\documentclass[final,3p,times]{elsarticle}
\usepackage{amssymb,amsmath,amsthm,epsfig,verbatim,url}
\usepackage[notref,notcite]{showkeys}
\usepackage{epstopdf}
\usepackage{mathrsfs}
\usepackage{ulem}
\usepackage{booktabs}
\biboptions{sort&compress}

\usepackage{tikz,pgfplots,tikz-3dplot}
\pgfplotsset{compat=1.7}
\usetikzlibrary{fit}

\usepackage[notref,notcite]{showkeys}
\newtheorem{thm}{Theorem}[section]
\newtheorem{lem}[thm]{Lemma}
\newtheorem{rem}[thm]{Remark}

\newcommand{\norm}[1]{\Vert#1\Vert}

\newcommand{\bR}{{\mathbb R}}
\newcommand{\bN}{{\mathbb N}}

\newcommand{\bI}{\mathbb{I}}
\newcommand{\bG}{\nabla}
\newcommand{\bD}{\mathbb{D}}
\newcommand{\bU}{\mathbb{U}}

\newcommand{\bP}{{\mathbb{P}}}

\newcommand{\mR}{\mathscr{R}}
\newcommand{\mS}{\mathcal{S}}
\newcommand{\drz}{{\rm d} r{\rm d}z}
\newcommand{\vol}{\operatorname{vol}}

\newcommand{\rd}{{\rm d}}

\newcommand{\id}{{\rm id}}
\newcommand{\dd}[1]{\frac{\rm d}{{\rm d}#1}}
\newcommand{\pp}[2]{\frac{\partial#2}{\partial#1}}
\newcommand{\ddt}{\dd{t}}
\newcommand{\nn}{\nonumber}
\newcommand{\ttau}{\Delta t}
\newcommand{\cira}{\mbox{$s\!\!\!\!\:/$}}
\newcommand{\mX}{\mathscr{X}}

\newcommand{\mZ}{{\mathfrak{x}}}
\newcommand{\Stab}{{\sf Stab}}
\newcommand{\Equi}{{\sf  Equi}}
\newcommand{\EquiV}{{\sf EquiV}}
\newcommand{\StabV}{{\sf  StabV}}

\makeatletter
\renewcommand*{\uuline}{%
  \bgroup
  \UL@setULdepth
  \markoverwith{%
    \lower\ULdepth\hbox{%
      \kern-.03em%
      \vtop{%
        \hrule width.2em%
        \kern 0.6pt % distance between the two underlines!
        \hrule
      }%
      \kern-.03em%
    }%
  }%
  \ULon
}
\makeatother
\newcommand{\mat}[1]{\smash{\uuline{#1}}}

\begin{document}
\begin{frontmatter}
\title{Arbitrary Lagrangian-Eulerian finite element approximations for axisymmetric two-phase flow}

\author[1]{Harald Garcke}
\address[1]{Fakult{\"a}t f{\"u}r Mathematik, Universit{\"a}t Regensburg, 
93040 Regensburg, Germany}
\ead{harald.garcke@ur.de}
\author[2]{Robert N\"urnberg}
\address[2]{Dipartimento di Mathematica, Universit\`a di Trento,
38123 Trento, Italy}
\ead{robert.nurnberg@unitn.it}
\author[3]{Quan Zhao}
\address[3]{School of Mathematical Sciences, University of Science and Technology of China, 230026 Hefei, China}
\ead{quanzhao@ustc.edu.cn}

\begin{abstract}
We analyze numerical approximations for axisymmetric two-phase flow in the arbitrary Lagrangian-Eulerian (ALE) framework. We consider a parametric formulation for the evolving fluid interface in terms of a one-dimensional generating curve. For the two-phase Navier-Stokes equations, we introduce both conservative and nonconservative ALE weak formulations in the 2d meridian half-plane. Piecewise linear parametric elements are employed for discretizing the moving interface, which is then coupled to a moving finite element approximation of the bulk equations. This leads to a variety of ALE methods, which enjoy either an equidistribution property or unconditional stability. Furthermore, we adapt these introduced methods with the help of suitable time-weighted discrete normals, so that  the volume of the two phases is exactly preserved on the discrete level.  Numerical results for rising bubbles and oscillating droplets are presented to show the efficiency and accuracy of these introduced methods. 
\end{abstract} 
\end{frontmatter}

\renewcommand{\thefootnote}{\arabic{footnote}}
\setcounter{equation}{0}
\setlength\parindent{24pt}

\section{Introduction}
Two-phase flows, and more generally multi-phase flows, occur in many natural phenomena and have wide applications in the oil and gas industries, engineering and scientific experiments. Numerical approximations of two-phase flow have been extensively studied in recent decades, and much effort has been devoted to accurate and efficient approximations of the evolving fluid interface. These include diffuse-interface methods \cite{Anderson1998, Grun2014two, Feng2006fully,Styles2008finite,GarckeHK16}, volume of fluid methods \cite{Hirt1981volume,Renardy02,Popinet09}, level set methods \cite{Sussman94level, Sethian99level, osher02level,olsson07}, and front-tracking methods \cite{UnverdiT92, Tryg01,Bansch01,Perot2003moving,Ganesan06,Quan07moving,BGN2013eliminating,BGN15stable,Agnese20,Zhao2021,Duan2022energy,GNZ23}.

Among these front-tracking approximations, one of the most prominent methods is the moving fitted mesh approach, where the discrete interface that separates the two fluids remains fitted to the bulk mesh. In particular, the interface mesh is approximated  by a lower-dimensional mesh and made up of faces of elements from the bulk mesh. This means that the bulk mesh needs to deform appropriately in time in order to match the evolving interface. The bulk equations are thus formulated in a moving frame of reference with a reference velocity which on the discrete level defines the movement of the bulk mesh. The natural way to employ this approach is the Lagrangian framework by simply prescribing the bulk mesh velocity according to the fluid velocity. However, this often leads to large distortions of the mesh due to the lack of control on the mesh\slash fluid velocity. A numerical method which allows for greater flexibility is the so-called arbitrary Lagrangian-Eulerian (ALE) method, where the reference velocity is somehow arbitrary and usually independent of the fluid velocity for the interior points. The original ALE approach was introduced in \cite{Noh64, Franck1964mixed,Hirt1974arbitrary} for hydrodynamic problems in the context of finite difference methods, and then generalized to free surface flows and fluid-structure interaction problems in the context of the finite element method, e.g., \cite{Belytschko78, Hughes81,donea82, Soulaimani1998}. The application of the ALE method to the two-phase flow can be found in e.g., \cite{Tryg01,Ganesan12arbitrary, Anjos20143d,Agnese20, Duan2022energy,GNZ23}. The main advantage of the ALE approach is the possibility to accurately capture the jumps of physical quantities across the interface, which on the contrary can be a major concern in the unfitted mesh approach \cite{Gross07extended, Frachon19cut, Claus2019cutfem}. This flexibility of the moving reference frame allows for excellent approximations in the case of small deformations. Nevertheless,  dynamic controls of the bulk and interface meshes are often necessary to prevent undesirable mesh distortions, especially when the interface exhibits strong deformations and topological changes. 

The stability of ALE finite element methods was first analyzed in \cite{Nobile1999} for the convection-diffusion equation based on either a conservative or nonconservative ALE formulation. In the paper, a condition of geometric conservative law (GCL) was proposed to guarantee an unconditional stability estimate that does not depend on the velocity of the ALE reference. A further stability analysis of ALE methods with a variety of time discretizations was considered in \cite{Boffi04stability}. We also refer the readers to \cite{Hron2006monolithic, Etienne09perspective, Liu2016second, Lozovskiy19, Lan2020, Wang2020energy, Kesler21} for the ALE methods of other problems such as the Stokes/parabolic moving interface problems and the fluid-structure intersection problem. As regards ALE methods for multi-phase flow, stability estimates can be found in e.g., \cite{Gerbeau2009generalized, ivanvcic2022energy, Duan2022energy, GNZ23}. For example, in \cite{ivanvcic2022energy},  the energy stability was established based on a GCL-type approximation of a conservative ALE formulation in the context of the single-phase of fluid. On assuming a divergence free velocity of the ALE frame, an energy-stable ALE method of the nonconservative form was recently proposed in \cite{Duan2022energy}. More recently in \cite{GNZ23}, the authors devised two structure-preserving ALE approximations for the two-phase incompressible flow in both the conservative and nonconservative form, and the introduced methods were shown to satisfy unconditional stability and exact volume preservation on the fully discrete level.

Despite the abundance of numerical work for two-phase flow, the computation for the fully 3d problem remains a very challenging task. The difficulties stem not only from the large size of the problem but also the complicated mesh manipulations which are often required in the front-tracking methods, see~\cite{Anjos20143d,BGN15stable, Duan2022energy}. Fortunately, in many situations the complex 3d problem can be reduced to a much simpler two-dimensional problem in the meridian half-plane provided that the considered flow satisfies rotational symmetry. In this axisymmetric setting the fluid interface can also be modeled by a one-dimensional generating curve, which dramatically reduces the computational complexity and troublesome work of mesh control. Existing numerical works for the axisymmetric two-phase flow can be found in Refs.~\cite{Sussman00coupled,Chessa2003,Kim05diffuse,Ganesan08accurate,Gros18interface,Huang22diffuse, GNZ23asy}. Very recently unfitted finite element approximations were analyzed by the authors in \cite{GNZ23asy}. In the current work, we aim to explore accurate and efficient numerical approximations for the axisymmetric two-phase flow in the ALE framework. In particular, the main results of this work are stated as follows. 
\begin{itemize}
\item Based on our recent 2d/3d work in \cite{GNZ23}, we introduce appropriate conservative and nonconservative ALE approximations which enable the stability of the fluid kinetic energy.
\item Inspired by the works in \cite{BGN19asy, BGN19variational}, we discuss two possible approximations of the surface tension forces which lead to either an unconditional stability estimate or to an equidistribution property.
\item Building on ideas in \cite{Jiang21,BGNZ22volume, GNZ23asy}, we adapt the introduced methods to further achieve volume-preserving approximations with the help of suitable time-integrated interface normals on the discrete level.
\end{itemize}

% organization of this paper

The remainder of the paper is organized as follows. In Section~\ref{sec:strongf} we discuss the axisymmetric setting for two-phase Navier-Stokes flow. Next, in Section~\ref{sec:ALEWF}, we introduce both the conservative and nonconservative ALE weak formulations for the axisymmetric flow and prove the volume preservation and energy stability for the weak solutions. A variety of finite element approximations of the introduced ALE formulations are then analyzed in Section~\ref{sec:FEM}. Subsequently, the solution method and numerical results for the introduced methods are presented in Section~\ref{sec:NR}. Finally we draw some conclusions in Section~\ref{sec:con}.

\section{The strong formulation}\label{sec:strongf}

\setcounter{equation}{0}

\begin{figure}[!htp]
\center
\definecolor{mcolor}{RGB}{127,0,255}
\newcommand{\AxisRotator}[1][rotate=0]{%
    \tikz [x=0.25cm,y=0.60cm,line width=.2ex,-stealth,#1] \draw (0,0) arc (-150:150:1 and 1);%
}
\begin{tikzpicture}[every plot/.append style={very thick}, scale = 1]
\begin{axis}[axis equal,axis line style=thick,axis lines=center, xtick style ={draw=none}, 
ytick style ={draw=none}, xticklabels = {}, 
yticklabels = {}, 
xmin=-0.1, xmax = 2.0, ymin = -0.4, ymax = 3.05]
after end axis/.code={  
   \node at (axis cs:0.0,2.5) {\AxisRotator[rotate=-90]};
   \draw[black, dashed, line width=1pt](axis cs: 1.8,0) --(axis cs:1.8,2.5);
   \draw[black, dashed, line width=1pt](axis cs: 0,2.5) --(axis cs:1.8,2.5);
   \draw[blue,->,line width=1.5pt] (axis cs:0,0) -- (axis cs:0.5,0);
   \draw[blue,->,line width=1.5pt] (axis cs:0,0) -- (axis cs:0,0.5);
   \node[blue] at (axis cs:0.5,-0.2){$\vec e_1$};
   \node[blue] at (axis cs:-0.2,0.5){$\vec e_2$};
   \draw[red,very thick] (axis cs: 0,0.6) arc[radius = 70, start angle= -90, end angle= 90];
   \draw[blue,->,line width=1.2pt] (axis cs:0.7*0.8660, 1.3+0.5*0.7) -- (axis cs:1.15*0.8660,1.3+0.5*1.15);
   \node[blue] at (axis cs:0.8,1.95){$\vec\nu$};
   \node[red] at (axis cs:0.8,0.85){$\Gamma(t)$};
   \node[mcolor] at (axis cs:0.3,1.2){$\mR_-(t)$};
   \node[mcolor] at (axis cs:1.0,1.4){$\mR_+(t)$};
   \node[black] at (axis cs:-0.2,1.4){$\partial_z\mR$};
}
\end{axis}
\end{tikzpicture} \qquad \qquad
\tdplotsetmaincoords{120}{50}
\begin{tikzpicture}[scale=2, tdplot_main_coords,axis/.style={->},thick]
\draw[axis] (-1, 0, 0) -- (1, 0, 0);
\draw[axis] (0, -1, 0) -- (0, 1, 0);
\draw[axis] (0, 0, -0.2) -- (0, 0, 2.7);
\draw[blue,->,line width=1.5pt] (0,0,0) -- (0,0.5,0) node [below] {$\vec e_1$};
\draw[blue,->,line width=1.5pt] (0,0,0) -- (0,0.0,0.5);
\draw[blue,->,line width=1.5pt] (0,0,0) -- (0.5,0.0,0);
\node[blue] at (0.2,0.4,0.1){$\vec e_3$};
\node[blue] at (0,-0.2,0.3){$\vec e_2$};
\node[red] at (0.9,0,1.4){$\mS(t)$};
\node at (0.0,0.0,2.3) {\AxisRotator[rotate=-90]};

\tdplottransformmainscreen{0}{0}{1.4}
\shade[tdplot_screen_coords, ball color = red] (\tdplotresx,\tdplotresy) circle (0.5);
\end{tikzpicture}
\caption{A sketch of the axisymmetric two-phase flow in the 2d meridian half-plane $\mR=\mR_-(t)\cup\mR_+(t)\cup\Gamma(t)$ as well as the unit vectors $\vec e_1$, $\vec e_2$ and $\vec e_3$, where $\Gamma(t)$ is the generating curve of the axisymmetric surface $\mS(t)$, and $\partial_z\mR$ is the artificial boundary of $\mR$ on $z$-axis.}
\label{fig:asy}
\end{figure}
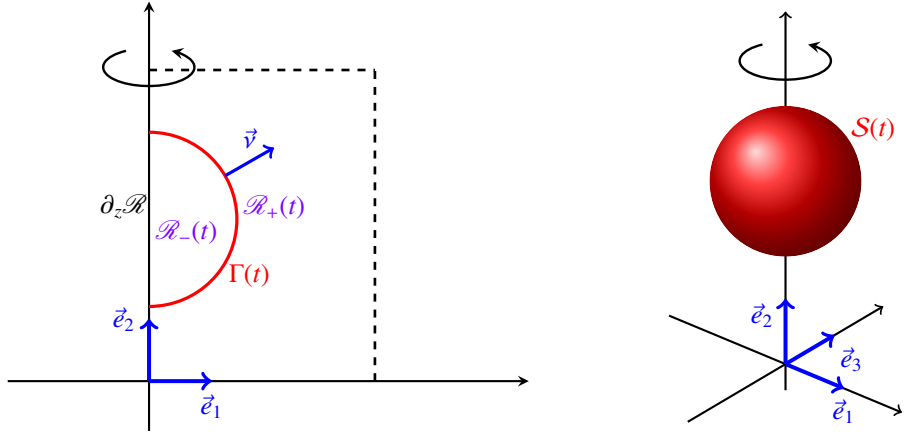

As shown in Fig.~\ref{fig:asy}, we assume that the two-phase flow satisfies rotational symmetry with respect to the $z$-axis. We then consider the problem in a bounded domain $\mR\subset\bR^2$ in the 2d meridian half-plane such that $\mR=\mR_+(t)\cup\mR_-(t)\cup\Gamma(t)$, where $\mR_\pm(t)$ correspond to the rotated sets for the domains of the two fluids, and $\Gamma(t)$ is the generating curve of the axisymmetric fluid interface $\mS(t)$. We further assume that there is no angular velocity and that the velocity components in the $r,z$-directions and the pressure are independent of the azimuthal angle. Thus we can introduce the velocity and pressure as variables in $\mR$
\begin{equation}
\vec u(\cdot, t)=\bigl(u^r(\cdot, t),\, u^z(\cdot,t )\bigr)^T: \mR\times[0,T]\to\bR^2\qquad\mbox{and}\qquad p(\cdot, t): \mR\times[0,T]\to\bR.\nn
\end{equation}
We denote by $\rho_\pm$ and $\mu_\pm$ the densities and viscosities of the two fluids, and introduce
\begin{equation}
\bG = \vec e_1\pp{r}{} + \vec e_2\pp{z}{},\qquad \mat{\bD}(\vec u)=\frac{1}{2}[\bG\vec u + (\bG\vec u)^T].\nn
\end{equation}
Then the axisymmetric Navier-Stokes equations in the two phases can be written as
\begin{subequations}\label{eqn:asyns}
\begin{alignat}{2}
\label{eq:asyns}
\rho_\pm\Bigl(\partial_t\vec u + [\vec u\cdot\bG]\vec u\Bigr)&=-\bG p + \frac{2}{r}\bG\cdot[r\,\mu_\pm\mat{\bD}(\vec u)]-\frac{2\mu_\pm (\vec u\cdot\vec e_1)\vec e_1}{r^2} + \rho_\pm \vec g\qquad &&\mbox{in}\quad\mR_\pm(t),\\
\frac{1}{r}\bG\cdot[r\vec u]& = \frac{1}{r}\pp{r}{(r\,u^r)}+ \pp{z}{u^z} =0\qquad &&\mbox{in}\quad\mR_\pm(t),
\label{eq:divfree}
\end{alignat}
\end{subequations}
where we denote $\vec g=(g^r,~g^z)^T$ as the body acceleration. Here \eqref{eq:asyns} and \eqref{eq:divfree} can be obtained via the axisymmetric reduction of the fully 3d governing equations for two-phase flow, where we refer the reader to e.g.,~\cite{bernardi99,GNZ23asy}.

We further assume that the generating curve $\Gamma(t)$ is an open curve with two end points which lie on the $z$-axis, so that the fluid interface $\mS(t)$ is a genus-0 surface without boundary. We introduce a parameterization of $\Gamma(t)$ as 
\begin{equation}
\vec\mZ(\cdot, t): \overline{\bI}\to\bR_{\geq 0}\times\bR\qquad\mbox{with}\quad\bI=(0,1),\quad \partial\bI=\{0,1\}.\nn
\end{equation}
The attachment of the two end points on the $z$-axis implies that 
\begin{subequations}\label{eqn:gammabd}
\begin{align}
\label{eq:gammabd1}
\vec\mZ(\alpha,t)\cdot\vec e_1=0\quad\forall\alpha\in\partial\bI,\\
\vec\mZ_\alpha(\alpha, t)\cdot\vec e_2 = 0\quad\forall\alpha\in\partial\bI,
\label{eq:gammabd2}
\end{align}
\end{subequations}
where \eqref{eq:gammabd2} is the axisymmetric condition meaning that the contact angle that $\Gamma(t)$ makes with the $z$-axis is $90^\circ$. We assume that $|\vec\mZ_\alpha(\cdot,t)|>0$ and $\vec \mZ(\cdot,t)\cdot\vec e_1>0$ in $\bI$, and the unit tangent and normal to the curve $\Gamma(t)$ are defined as 
\begin{equation}\label{eq:taunu}
\vec\tau(\alpha,t)  = \vec\mZ_s = |\vec\mZ_\alpha|^{-1}\vec\mZ_\alpha,\qquad \vec\nu = -(\vec\tau)^\perp,
\end{equation}
where $s$ is the arc length of $\Gamma(t)$ with  $\partial_s =
|\vec{\mZ}_\alpha|^{-1}\,\partial_\alpha$, and $(\cdot)^\perp$ denotes a clockwise rotation of a vector by $\frac{\pi}{2}$. Besides, we introduce the mean curvature of $\mS(t)$
\begin{equation}\label{eq:asykappa}
\varkappa = \kappa - (\vec\mZ\cdot\vec e_1)^{-1}(\vec\nu\cdot\vec e_1)\qquad\mbox{with}\quad\kappa\vec\nu = \vec\mZ_{ss},
\end{equation}
where $\kappa$ is the curvature of the generating curve $\Gamma(t)$.

We have interface conditions on $\mS(t)$, which lead to 
\begin{equation}\label{eq:ifc}
[\vec u]_-^+ = \vec 0,\qquad [2\,\mu_\pm\mat{\bD}(\vec u) - p\mat{Id}]_-^+\,\vec\nu = -\gamma\varkappa\,\vec\nu,\qquad \partial_t\vec\mZ\cdot\vec\nu = \vec u\cdot\vec\nu\qquad\mbox{on}\quad\Gamma(t),
\end{equation}
where $[\cdot]_-^+$ denotes the jump values from $\mR_-(t)$ to $\mR_+(t)$, $\gamma$ is the surface tension of the fluid interface, and $\mat{Id}\in\bR^{2\times 2}$ is the identity matrix. 

Furthermore, the boundary of $\mR$ is given by $\partial\mR=\partial\mR_1\cup\partial\mR_2\cup\partial_z\mR$, where $\partial_z\mR$ is the artificial boundary of $\mR$ on the $z$-axis. We impose a no-slip condition on $\partial_1\mR$ and a free-slip condition on $\partial_2\mR$ as
\begin{subequations}\label{eqn:Rbd}
\begin{alignat}{2}
\label{eq:noslip}
\vec u &= \vec 0\qquad &&\mbox{on}\quad\partial_1\mR,\\
\vec u\cdot\vec n = 0,\quad\mat{\bD}(\vec u)\,\vec n\cdot\vec t &= 0\qquad &&\mbox{on}\quad\partial_2\mR,\label{eq:freeslip}
\end{alignat}
\end{subequations}
where $\vec n$ and $\vec t$ are normal and tangent to $\partial_2\mR$. On the artificial boundary $\partial_z\mR$, we need
\begin{equation}\label{eq:ARbd}
\vec u \cdot\vec e_1 = 0,\qquad \pp{r}{\vec u}\cdot\vec e_2 = 0\qquad\mbox{on}\quad\partial_z\mR.
\end{equation}
Here the first condition in \eqref{eq:ARbd} is enforced on recalling the term $-\frac{2\mu_\pm(\vec u\cdot\vec e_1)\vec e_1}{r^2}$ in \eqref{eq:asyns}, while the second condition results from the axisymmetry. It is not difficult to show that \eqref{eq:ARbd} is equivalent to the free-slip condition.

\section{The ALE weak formulations}\label{sec:ALEWF}
We introduce a family of  ALE mappings $\{\mathcal{\vec A}[t]\}_{t\in[0,T]}$ such that
\begin{equation}
\vec{\mathcal{A}}[t]: \mathcal{O}\to \mR, \qquad \vec y \mapsto \vec{\mathcal{A}}[t](\vec y) = \vec x(\vec y, t)\quad \mbox{for all}\quad t\in[0,T],\quad \vec y\in\mathcal{O}.
\end{equation}
We further assume that $\vec{\mathcal{A}}[t]\in [W^{1,\infty}(\mathcal{O})]^2$ and $\vec{\mathcal{A}}[t]^{-1}\in [W^{1,\infty}(\mR)]^2$, and define the mesh velocity 
\begin{equation}
\vec w(\vec x, t) :=\left.\frac{\partial\vec x(\vec y,t)}{\partial t} \right|_{\vec y = \vec{\mathcal{A}}[t]^{-1}(\vec x)}\quad\mbox{for all}\quad t\in[0,T],\quad \vec x\in\mR.\label{eq:ALEmeshv}
\end{equation}
On the interface and boundary, the mesh velocity is required to satisfy
\begin{equation}
(\vec w-\vec u)\cdot\vec\nu = 0\quad\mbox{on}\quad\Gamma(t),\qquad(\vec w-\vec u)\cdot\vec n = 0\quad\mbox{on}\quad\partial\mR. 
\label{eq:meshvcond}
\end{equation}
For a vector field $\vec\varphi:\mR\times[0,T]\to\bR^2$, we also introduce the derivative with respect to the moving ALE frame as
\begin{equation}
\partial_t^\circ\vec\varphi = \partial_t\vec\varphi + [\vec w\cdot\bG]\vec\varphi.\label{eq:meshD}
\end{equation}

We define $(\cdot,\cdot)$ as the $L^2$-inner product on $\mR$ and denote by
\begin{align}
L_a^2(\mR) &:= \{\chi: (r^a,~\,\chi^2) < +\infty\}\quad a\in\bN,\nn\\
H_a^1(\mR)&:=\{\chi: \chi\in L_a^2(\mR),\;\bG\chi\in [L_a^2(\mR)]^2\}\nn 
\end{align}
the weighted $L^2$ and $H^1$ spaces over $\mR$, respectively.  We then introduce the function spaces for the velocity and pressure as
\begin{subequations}
\begin{align}
\bU & := \Bigl\{\vec\chi\in[H_1^1(\mR)]^2\;:\;(\vec\chi\cdot\vec e_1)\in L^2_{-1}(\mR),\; \vec\chi=\vec 0\;\mbox{on}\;\partial_1\mR,\;\vec\chi\cdot\vec n = 0\;\mbox{on}\;\partial_2\mR\Bigr\},\\
\mathbb{V}&:=H^1(0,T;\,[L^2_1(\mR)]^2)\cap L^2(0,T;\bU),\qquad\bP := \bigl\{\chi\in L^2_1(\mR):  (r,\,\chi)=0\bigr\}.
\end{align}
\end{subequations}
Besides, we denote by $\langle\cdot,\cdot\rangle$ the $L^2$-inner product on $\bI$, and define the function space 
\begin{equation}
V_{\partial} = \left\{\vec\eta\in [H^1(\bI)]^2\;:\;\vec\eta\cdot\vec e_1 = 0\;\mbox{on}\;\partial\bI\right\}.\nn
\end{equation}

\subsection{Nonconservative ALE formulations}
Denote
\begin{equation}
\rho(\cdot,t) = \rho_+\mX_{\mR_+(t)} + \rho_-\mX_{\mR_-(t)},\qquad \mu(\cdot,t)  = \mu_+\mX_{\mR_+(t)} + \mu_-\mX_{\mR_-(t)},
\end{equation}
where $\mX_E$ is the characteristic function of the set $E$. For the pressure and viscous term in \eqref{eq:asyns}, we take the inner product with $\vec\chi\,r$ for $\vec\chi\in\bU$. This leads to
\begin{align}
&\bigl(-\bG p,~\vec\chi\, r\bigr) + \bigl(\bG\cdot[2r\mu\mat{\bD}(\vec u)],~\vec\chi\bigr) -2\bigl(\mu\,r^{-1}\,[\vec u\cdot\vec e_1],~[\vec\chi\cdot\vec e_1]\bigr)\nn\\
&= \bigl(p,~\bG[r\,\vec\chi]\bigr)-2\bigl(\mu\,\mat{\bD}(\vec u),~\mat{\bD}(\vec\chi)\,r\bigr)-2\bigl(\mu\,r^{-1}\,[\vec u\cdot\vec e_1],~[\vec\chi\cdot\vec e_1]\bigr)+\gamma\int_{\Gamma(t)}(\vec\mZ\cdot\vec e_1)\varkappa\vec\nu\cdot\vec\chi\,\rd s, 
\label{eq:ALE1}
\end{align}
where we used integration by parts, the interface condition in \eqref{eq:ifc}, and the boundary conditions \eqref{eqn:Rbd} and \eqref{eq:ARbd}.

Denote by $(\cdot,\cdot)_{\mR_\pm(t)}$ the $L^2$-inner products over $\mR_\pm(t)$, respectively. For the inertia term in \eqref{eq:asyns}, on recalling \eqref{eq:meshD} we have that
\begin{equation}\label{eq:ALE2}
\bigl(\partial_t\vec u + [\vec u \cdot\bG]\vec u, ~\vec\chi\,r\bigr)_{\mR_\pm(t)}= \bigl(\partial_t^\circ\vec u, ~\vec\chi\,r\bigr)_{\mR_\pm(t)} + \bigl([\vec u-\vec w] \cdot\bG\vec u,~\vec\chi\,r\bigr)_{\mR_\pm(t)}\quad\forall\vec\chi\in[H^1(\mR)]^2.
\end{equation}
We can rewrite 
\begin{align}
\bigl([\vec u-\vec w] \cdot\bG\vec u,~\vec\chi\,r\bigr)_{\mR_\pm(t)}&=\frac{1}{2}\left[\bigl(([\vec u-\vec w]\cdot\bG)\vec u,~\vec \chi\,r\bigr)_{\mR_\pm(t)}-\bigl(([\vec u-\vec w]\cdot\bG)\vec\chi,~\vec u\,r\bigr)_{\mR_\pm(t)}\right]+\frac{1}{2}\bigl(\vec u - \vec w,~\bG(\vec u\cdot\vec\chi)r\bigr)_{\mR_\pm(t)}\nn\\
&=\frac{1}{2}\left[\bigl(([\vec u-\vec w]\cdot\bG)\vec u,~\vec \chi\,r\bigr)_{\mR_\pm(t)}-\bigl(([\vec u-\vec w]\cdot\bG)\vec\chi,~\vec u\,r\bigr)_{\mR_\pm(t)}\right] + \frac{1}{2}\bigl(\bG\cdot[r\,\vec w],~\vec u\cdot\vec \chi\bigr)_{\mR_\pm(t)},\label{eq:ALE3n}
\end{align}
where the last equality results from the divergence free condition in \eqref{eq:divfree} and the boundary condition in \eqref{eq:meshvcond}.
We then multiply \eqref{eq:ALE1} with $\rho_\pm$ and combine these two equations. This yields that
\begin{equation}
\bigl(\rho\,(\partial_t\vec u +[\vec u\cdot\bG]\vec u),~\vec\chi\,r\bigr)=\bigl(\rho\,\partial_t^\circ\vec u,~\vec\chi\,r\bigr) + \frac{1}{2}\bigl(\rho\,\bG\cdot[r\,\vec w],~\vec\chi\bigr) + \mathscr{A}(\rho,\vec u-\vec w;\vec u,\vec\chi),\label{eq:ALE3}
\end{equation}
where $\mathscr{A}(\rho,\vec v;\vec u,\vec\chi)$ is the antisymmetric term 
\begin{equation}
\mathscr{A}(\rho,\vec v;\vec u,\vec\chi) = \frac{1}{2}\left[\bigl(\rho\,(\vec v\cdot\nabla)\vec u,~\vec \chi\,r\bigr)-\bigl(\rho\,(\vec v\cdot\nabla)\vec \chi,~\vec u\,r\bigr)\right].\nn
\end{equation}

It was shown in \cite{BGN19asy} that taking the inner product of \eqref{eq:asykappa} with $(\vec\mZ\cdot\vec e_1)\vec\nu\cdot\vec\eta\,|\vec\mZ_\alpha|$ on $\bI$ with $\vec\eta\in V_{\partial}$ leads to
\begin{equation} \label{eq:ALE4}
\big\langle(\vec\mZ\cdot\vec e_1)\,
\varkappa\,\vec\nu,~\vec\eta\, |\vec x_\alpha|\big\rangle 
+\big\langle\vec\eta \cdot\vec e_1,~|\vec\mZ_\alpha|\big\rangle 
+ \big\langle(\vec\mZ\cdot\vec e_1)\vec \mZ_\alpha,~\vec\eta_\alpha|\vec x_\alpha|^{-1}\big\rangle  
= 0\qquad \forall\ \vec\eta\in V_{\partial}\,.
\end{equation}
Collecting the results in \eqref{eq:ALE1}, \eqref{eq:ALE3} and \eqref{eq:ALE4}, we then introduce the following  ALE weak formulation in the nonconservative form. Let the initial velocity $\vec u_0\in\bU$ and interface parameterization $\vec\mZ_0\in V_{\partial}$. For $t\in(0,T]$, we find $\vec u(\cdot, t)\in\mathbb{V}$, $p\in\bP$, $\vec\mZ(\cdot,t)\in\ V_{\partial}$ and $\varkappa(\cdot,t)\in L^2(\bI)$ such that 
\begin{subequations}\label{eqn:weakform}
\begin{align}
&\bigl(\rho\partial_t^\circ\vec u,~\vec\chi\,r\bigr) + \frac{1}{2}\bigl(\rho\bG\cdot[r\vec w], ~\vec u\cdot\vec\chi\bigr)+\mathscr{A}(\rho,\vec u-\vec w; \vec u, \vec\chi)  + 2\bigl(\mu\,r^{-1}\,[\vec u\cdot\vec e_1],~[\vec\chi\cdot\vec e_1]\bigr)- \bigl(p,~\bG\cdot[r\,\vec\chi]\bigr)\nn\\
&\hspace{0.5cm}+ 2\Bigl(\mu\,\mat{\bD}(\vec u),~\mat{\bD}(\vec\chi)\,r\Bigr) -\,\gamma\big\langle\,(\vec\mZ\cdot\vec e_1)\,\varkappa,~\vec\nu\cdot\vec\chi\,|\vec\mZ_\alpha|\,\big\rangle = \bigl(\rho\,r\,\vec g,~\vec\chi\bigr)\qquad\forall\vec\chi\in\mathbb{V},\label{eq:weak1}\\[0.3em]
\label{eq:weak2}
&\hspace{0.2cm}\bigl(\bG\cdot[r\,\vec u], ~q\bigr)=0\qquad\forall q\in \bP,\\[0.3em]
\label{eq:weak3}
&\hspace{0.2cm}\big\langle(\vec\mZ\cdot\vec e_1)\,\vec \mZ_t\cdot\vec\nu,~\zeta\,|\vec \mZ_\alpha|\big\rangle-\big\langle(\vec\mZ\cdot\vec e_1)\,\vec u\cdot\vec\nu,~\zeta\,|\vec \mZ_\alpha|\big\rangle = 0\qquad\forall \zeta\in L^2(\bI),\\[0.3em]
\label{eq:weak4}
&\hspace{0.2cm}\big\langle(\vec\mZ\cdot\vec e_1)\,\varkappa\,\vec\nu,~\vec\eta\,|\vec\mZ_\alpha|\big\rangle+\big\langle\vec\eta\cdot\vec e_1,~|\vec\mZ_\alpha|\big\rangle + \big\langle(\vec\mZ\cdot\vec e_1)\,\vec \mZ_\alpha,~\vec\eta_\alpha\,|\vec\mZ_\alpha|^{-1}\big\rangle=0\qquad\forall\vec\eta\in\ V_{\partial}.
\end{align}
\end{subequations}
We note that \eqref{eq:weak1} results from \eqref{eq:ALE1} and \eqref{eq:ALE2},  \eqref{eq:weak2} and \eqref{eq:weak3} are due to the incompressibility condition \eqref{eq:divfree} and the kinetic equation in \eqref{eq:ifc}, while \eqref{eq:weak4} is a direct result of \eqref{eq:ALE4}.

For simplicity we denote $\vec\mZ(\cdot, t)= \vec\mZ(t)$ and $\vec u(\cdot, t) = \vec u(t)$.  We then introduce $A(\vec\mZ(t))$ and $\vol(\vec\mZ(t))$ as the surface area and the enclosed volume of $\mS(t)$, respectively:
\begin{subequations}\label{eqn:SAV}
\begin{align}\label{eq:area}
A(\vec\mZ(t))&=\int_{\mS(t)}1\rd A = 2\pi\int_\bI(\vec\mZ\cdot\vec e_1)|\vec\mZ_\alpha|\,\rd\alpha,\\
\vol(\vec\mZ(t)) &=\pi\int_{\Gamma(t)}(\vec\mZ\cdot\vec e_1)^2(\vec\nu\cdot\vec e_1)\rd s= \pi\int_{\bI}(\vec\mZ\cdot\vec e_1)^2(\vec\nu\cdot\vec e_1)\,|\vec\mZ_\alpha|\,\rd\alpha,\label{eq:vol}
\end{align}
\end{subequations}
where for \eqref{eq:vol} the reader can refer to \cite[(3.10)]{BGN19asy}. The total free energy of the system is given by
\begin{equation}
\label{eq:Energy}
E(\rho, \vec u(t),\vec\mZ(t))=\pi \int_{\mR}\rho\,|\vec u|^2\,r\,\drz + \gamma A(\vec\mZ(t)).
\end{equation}
It is straightforward to obtain the time derivative of \eqref{eq:area} as
\begin{subequations}
\begin{align}
\label{eq:ddtA}
\ddt A(\vec\mZ(t)) &=  2\pi\int_{\bI}\left[(\vec\mZ_t\cdot\vec e_1)\,|\vec\mZ_\alpha| + (\vec\mZ\cdot\vec e_1)\,(\vec\mZ_t)_\alpha\cdot\vec\mZ_\alpha|\vec\mZ_\alpha|^{-1}\right]\rd\alpha.
\end{align}
Moreover, it follows form the Reynolds transport theorem that (see also \cite[(2.22)]{BGN19asy})
\begin{align}
\ddt \vol(\vec\mZ(t)) =2\pi\int_{\Gamma(t)}(\vec\mZ\cdot\vec e_1)\,(\vec\mZ_t\cdot\vec\nu)\,\rd s=2\pi\int_{\bI}(\vec\mZ\cdot\vec e_1)(\vec\mZ_t\cdot\vec\nu)|\vec\mZ_\alpha|\,\rd\alpha.
\label{eq:ddtM}
\end{align}
\end{subequations}

 Now choosing $\vec\chi = \vec u$ in \eqref{eq:weak1}, $q = p$ in \eqref{eq:weak2}, $\zeta = \gamma\varkappa$ in \eqref{eq:weak3} and $\vec\eta = \vec\mZ_t$ in \eqref{eq:weak4} and combing these equations yields that
\begin{align}
\label{eq:energylaw}
\frac{1}{2\pi}\ddt E(\rho,\vec u(t),\vec\mZ(t)) &= -2\Bigl(\mu\,\mat{\bD}(\vec u),\mat{\bD}(\vec u)\,r\Bigr) + 2\Bigl(r^{-1}\,\mu,~[\vec u\cdot\vec e_1]^2\Bigr) + \Bigl(\rho\,r\,\vec u,~\vec g\Bigr),
\end{align}
on recalling \eqref{eq:ddtA} and \eqref{eq:ddtK}. Besides,  it is not difficult to show that
 \begin{equation}
 \mX_{_{\mR_-(t)}}-\omega(t)\in\bP\quad\mbox{with}\quad \omega(t) = \frac{\int_{\mR_-(t)}r\drz}{\int_\mR r\drz},\quad t\in[0,T].\nn
 \end{equation}
Then setting $q = \mX_{_{\mR_-(t)}}-\omega(t)\in\bP$ in \eqref{eq:weak2} and $\zeta=1$ in \eqref{eq:weak3}, we arrive at
\begin{equation}
\ddt \vol(\vec\mZ(t))=2\pi\big\langle(\vec\mZ\cdot\vec e_1)(\vec\mZ_t\cdot\vec\nu),~|\vec\mZ_\alpha|\big\rangle = 2\pi\big\langle(\vec\mZ\cdot\vec e_1)\vec u\cdot\vec\nu,~|\vec\mZ_\alpha|\big\rangle=2\pi\Bigl(\nabla\cdot[r\,\vec u],~1\Bigr)=0. \label{eq:vollaw}
\end{equation}
This implies the energy stability and volume conservation within the weak formulation.

\medskip
On recalling \eqref{eq:asykappa}, it is also possible to consider an ALE weak formulation which treats the curvature of $\Gamma(t)$ as an unknown. Let $\vec u_0\in\bU$ and $\vec\mZ_0\in V_{\partial}$. For $t\in(0,T]$, we find $\vec u(\cdot, t)\in\mathbb{V}$, $p\in\bP$, $\vec\mZ(\cdot,t)\in\ V_{\partial}$ and $\kappa(\cdot,t)\in L^2(\bI)$ such that 
\begin{subequations}\label{eqn:Aweakform}
\begin{align}
&\bigl(\rho\partial_t^\circ\vec u,~\vec\chi\,r\bigr) + \frac{1}{2}\bigl(\rho\bG\cdot[r\vec w], ~\vec u\cdot\vec\chi\bigr)+\mathscr{A}(\rho,\vec u-\vec w; \vec u, \vec\chi)  + 2\bigl(\mu\,r^{-1}\,[\vec u\cdot\vec e_1],~[\vec\chi\cdot\vec e_1]\bigr)- \bigl(p,~\bG\cdot[r\,\vec\chi]\bigr)\nn\\
&\hspace{0.5cm}+ 2\Bigl(\mu\,\mat{\bD}(\vec u),~\mat{\bD}(\vec\chi)\,r\Bigr) -\,\gamma\big\langle\,(\vec\mZ\cdot\vec e_1)\,\kappa - \vec\nu\cdot\vec e_1,~\vec\nu\cdot\vec\chi\,|\vec\mZ_\alpha|\,\big\rangle = \bigl(\rho\,r\,\vec g,~\vec\chi\bigr)\qquad\forall\vec\chi\in\mathbb{V},\label{eq:Aweak1}\\[0.3em]
\label{eq:Aweak2}
&\hspace{0.2cm}\bigl(\bG\cdot[r\,\vec u], ~q\bigr)=0\qquad\forall q\in \bP,\\[0.3em]
\label{eq:Aweak3}
&\hspace{0.2cm}\big\langle(\vec\mZ\cdot\vec e_1)\,\vec \mZ_t\cdot\vec\nu,~\zeta\,|\vec \mZ_\alpha|\big\rangle-\big\langle(\vec\mZ\cdot\vec e_1)\,\vec u\cdot\vec\nu,~\zeta\,|\vec \mZ_\alpha|\big\rangle = 0\qquad\forall \zeta\in L^2(\bI),\\[0.3em]
\label{eq:Aweak4}
&\hspace{0.2cm}\big\langle\kappa\,\vec\nu,~\vec\eta\,|\vec\mZ_\alpha|\big\rangle + \big\langle\,\vec \mZ_\alpha,~\vec\eta_\alpha\,|\vec\mZ_\alpha|^{-1}\big\rangle=0\qquad\forall\vec\eta\in\ V_{\partial}.
\end{align}
\end{subequations}
Similarly, choosing $\vec\chi= \vec u$ in \eqref{eq:Aweak1}, $q=p$ in \eqref{eq:Aweak2}, $\zeta= \gamma[\kappa - (\vec\mZ\cdot\vec e_1)^{-1}(\vec\nu\cdot\vec e_1)]$ and $\vec\eta=(\vec\mZ\cdot\vec e_1)\vec\mZ_t$ in \eqref{eq:Aweak4} leads to \eqref{eq:energylaw} as well. Besides, \eqref{eq:vollaw} is straightforward by choosing $q = \mX_{_{\mR_-(t)}}-\omega(t)$ in \eqref{eq:weak2} and $\zeta=1$ in \eqref{eq:weak3}

In fact, \eqref{eqn:weakform} can also be obtained by using an axisymmetric reduction of the 3d ALE formulation in \cite[(4.9)]{GNZ23}. Compared to \eqref{eqn:Aweakform} the formulation \eqref{eqn:weakform} has the advantage that on discretization it leads to an unconditionally stable method. On the other hand, the discretization of \eqref{eqn:Aweakform}, thanks to \eqref{eq:Aweak4}, will admit an equidistribution property, similar to that in the 2d case \cite{BGN07}.

\subsection{Conservative ALE formulations}

On recalling \eqref{eq:ddtchiu} and \eqref{eq:ALE3}, it is not difficult to obtain 
\begin{equation}
\Bigl(\rho\,\partial_t^\circ\vec u,~\vec\chi\,r\Bigr) + \frac{1}{2}\Bigl(\rho\nabla\cdot[r\vec w],~\vec u\cdot\vec\chi\Bigr)=\ddt\Bigl(\rho\vec u,~\vec\chi\,r\Bigr) - \frac{1}{2}\Bigl(\rho\nabla\cdot[r\vec w],~\vec u\cdot\vec\chi\Bigr) - \Bigl(\rho\partial_t^\circ\vec\chi,~\vec u\,r\Bigr)\quad\forall\vec\chi\in[H^1(\mR)]^2.\label{eq:ncr}
\end{equation}
Then we could introduce an ALE weak formulation in the conservative form by modifying \eqref{eqn:weakform} slightly.  Precisely, we replace \eqref{eq:weak1} with 
\begin{align}
&\ddt\Bigl(\rho\vec u,~\vec\chi\,r\Bigr) - \frac{1}{2}\Bigl(\rho\nabla\cdot[r\vec w],~\vec u\cdot\vec\chi\Bigr) - \Bigl(\rho\partial_t^\circ\vec\chi,~\vec u\,r\Bigr) +\mathscr{A}(\rho,\vec u-\vec w; \vec u, \vec\chi)  + 2\bigl(\mu\,r^{-1}\,[\vec u\cdot\vec e_1],~[\vec\chi\cdot\vec e_1]\bigr)\nn\\
&\hspace{0.5cm}- \bigl(p,~\bG\cdot[r\,\vec\chi]\bigr)+ 2\bigl(\mu\,\mat{\bD}(\vec u),~\mat{\bD}(\vec\chi)\,r\bigr) -\,\gamma\big\langle\,(\vec\mZ\cdot\vec e_1)\,\varkappa,~\vec\nu\cdot\vec\chi\,|\vec\mZ_\alpha|\,\big\rangle = \bigl(\rho\,r\,\vec g,~\vec\chi\bigr)\qquad\forall\vec\chi\in\mathbb{V}.\label{eq:weakc}
\end{align}
Similarly for \eqref{eqn:Aweakform}, we could obtain a new formulation in the conservative form by replacing \eqref{eq:Aweak1} with 
\begin{align}
&\ddt\Bigl(\rho\vec u,~\vec\chi\,r\Bigr) - \frac{1}{2}\Bigl(\rho\nabla\cdot[r\vec w],~\vec u\cdot\vec\chi\Bigr) - \Bigl(\rho\partial_t^\circ\vec\chi,~\vec u\,r\Bigr) +\mathscr{A}(\rho,\vec u-\vec w; \vec u, \vec\chi)  + 2\bigl(\mu\,r^{-1}\,[\vec u\cdot\vec e_1],~[\vec\chi\cdot\vec e_1]\bigr)\nn\\
&\hspace{0.5cm}- \bigl(p,~\bG\cdot[r\,\vec\chi]\bigr)+ 2\bigl(\mu\,\mat{\bD}(\vec u),~\mat{\bD}(\vec\chi)\,r\bigr) -\,\gamma\big\langle\,(\vec\mZ\cdot\vec e_1)\,\kappa - (\vec\nu\cdot\vec e_1),~\vec\nu\cdot\vec\chi\,|\vec\mZ_\alpha|\,\big\rangle = \bigl(\rho\,r\,\vec g,~\vec\chi\bigr)\qquad\forall\vec\chi\in\mathbb{V}.\label{eq:Aweakc}
\end{align}

 In a similar manner to \eqref{eqn:weakform} and \eqref{eqn:Aweakform}, it is not difficult to prove the energy law \eqref{eq:energylaw} for the two conservative ALE weak formulations in view of \eqref{eq:ncr}. In fact,  \eqref{eq:weakc} can also be derived via an axisymmetric reduction of the 3d conservative ALE formulation in \cite[(4.12)]{GNZ23}.

\section{Finite element approximations}\label{sec:FEM}
\setcounter{equation}{0}

In this section, we propose ALE finite element approximations for the introduced four weak formulations in \S\ref{sec:ALEWF}, and explore the properties of these approximating methods.

\subsection{The discretization}
We employ a uniform partition of the time interval and the reference domain $\bI$ as follows
\begin{align}
&[0,T]=\cup_{m=1}^M[t_{m-1},t_{m}]\quad\mbox{with}\quad t_m = m\ttau,\quad \ttau = T/M,\nn\\
&\bI =\cup_{j=1}^{J_\Gamma}[\alpha_{j-1},\alpha_{j}]\quad\mbox{with}\quad \alpha_j = j\,h,\quad h = 1/J_{\Gamma}.\nn
\end{align}
We define the finite element spaces on $\bI$ as
\[
V^h := \bigl\{\zeta \in C(\overline {\bI}) : \zeta\!\mid_{\bI_j} \nn\
\text{is affine}\ \forall\ j=1,\ldots, J_\Gamma\bigr\},\qquad V^h_{\partial}=V_\partial\cup[V^h]^2.
\]
Let $\{\vec\mZ^m\}_{0\leq m\leq M}$ be an approximation to $\{\vec\mZ(t)\}_{t\in[0,T]}$ with $\vec\mZ^m\in V_\partial^h$.  We define the polygonal curve $\Gamma^m = \vec X^m(\overline{\bI})$. Throughout this section we assume
that
\begin{equation*} 
\vec X^m \cdot\vec e_1 > 0 \quad \text{in }\
\overline {\bI}\setminus \partial \bI
\quad\text{and}\quad
|\vec{X}^m_\alpha| > 0 \quad \text{in } \bI
\qquad 0\leq m\leq M,
\end{equation*}
so that we can set 
\begin{equation*} 
\vec\tau^m = \vec X^m_s = \frac{\vec X^m_\alpha}{|\vec X^m_\alpha|} 
\qquad \mbox{and} \qquad \vec\nu^m = -(\vec\tau^m)^\perp\,.
\end{equation*}

Let 
\begin{equation}
\overline{\mR}=\cup_{o\in\mathscr{T}^m}\overline{o}\quad\mbox{with}\quad\mathscr{T}^m=\{o_j^m: j = 1,\cdots, J_\mR\},\qquad\mbox{and}\quad Q^m=\{\vec q_k^m: k = 1,\cdots, K\}, \nn
\end{equation}
be a regular partition at time $t=t_m$, where $\mathscr{T}^m$ are mutually disjoint and non-degenerate triangles, and $Q^m$ is the set of vertices of $\mathscr{T}^m$. We introduce the finite element spaces associated with $\mathscr{T}^m$ as
\begin{subequations}
\begin{align}
S_k^m&:=\left\{\varphi\in C(\overline{\mR}):\; 
\varphi|_{o_j^m}\in \mathcal{P}_k(o_j^m),\;\forall\ j=1,\cdots,J_{\mR}\right\},\quad k\in\bN_+, \\
S_0^m&:=\{\varphi\in L^2(\Omega):\; \varphi|_{o_j^m}\;\mbox{is constant},\; \forall j=1,\cdots,J_{\mR}\},
\end{align}
\end{subequations}
where $\mathcal{P}_k(o_j^m)$ denotes the space of polynomials of degree $k$ on $o_j^m$.

We use the moving fitted finite element method so that the discretization of the interface $\Gamma^m$ is fitted to the triangulation of $\mR$,
i.e. the line segments making up $\Gamma^m$ are all edges of elements in $\mathscr{T}^m$. Let $\mR_-^m$ be the interior region enclosed by $\Gamma^m$ and $\partial_z\mR$, and $\mR_+^m$ be the exterior region. We then approximate $\mu(\cdot,t)$ and $\rho(\cdot, t)$ with $\rho^m\in S_0^m$ and $\mu^m\in S_0^m$ such that
\begin{equation}\label{eq:Ddensity}
\rho^m = \rho_-\mX_{_{\mR_-^m}} + \rho_+\mX_{_{\mR_+^m}},\qquad \mu^m = \mu_-\mX_{_{\mR_-^m}} + \mu_+\mX_{_{\mR_+^m}}.\nonumber
\end{equation}
We denote by $\bU^m$ and $\bP^m$ the velocity and pressure approximation spaces, respectively. In the present work we consider the following pair element \cite{BGN2013eliminating, BGN15stable, GNZ23}
\begin{equation}
\label{eq:P2P1P0}
\mbox{P2-(P1+P0)}:\quad\bigl(\mathbb{U}^m,~\mathbb{P}^m\bigr)= \bigl([S_2^m]^2\cap\mathbb{U},~(S_1^m+S_0^m)\cap\mathbb{P}\bigr),
\end{equation}
which is able to capture the pressure jump across the interface.  

\subsection{Discrete ALE mappings}
\label{sec:dale}

In order to match the evolving polygonal curves $\Gamma^m$, the bulk mesh $\mathscr{T}^m$ needs to be constructed appropriately. For each $m\geq 1$, we now assume that we are given the polyhedral curve $\Gamma^{m}=\vec X^m(\bI)$. We then construct $\mathscr{T}^{m}$ based on $\mathscr{T}^{m-1}$ by simply moving the vertices in $Q^{m-1}$ according to the displacement vectors
\begin{equation}\label{eq:meshQU}
\vec q_k^{m} = \vec q_k^{m-1} + \vec\psi(\vec q_k^{m-1}),\qquad 1\leq k\leq K, \quad 1\leq m\leq M,
\end{equation}
while preserving the connectivity and topology of the bulk mesh.  On introducing 
\begin{align}
\mathbb{Y}^{m-1}&=\bigl\{\vec\chi\in [S_1^{m-1}]^2,\,\vec\chi\cdot\vec n = 0\;\;\mbox{on}\;\;\partial\mR;\;\vec\chi = \vec X^m - \vec X^{m-1}\;\;\mbox{on}\;\Gamma^{m-1}\bigr\},\nn\\
\mathbb{Y}_0^{m-1}&=\bigl\{\vec\chi\in [S_1^{m-1}]^2,\,\vec\chi\cdot\vec n = 0\;\;\mbox{on}\;\;\partial\mR;\;\vec\chi = \vec 0\;\;\mbox{on}\;\;\Gamma^{m-1}\bigr\},\nn
\end{align}
we solve for the displacement $\vec\psi^m\in\mathbb{Y}^{m-1}$ via an elastic equation such that 
\begin{equation}
2\bigl(\lambda^{m-1}\,\mat{\bD}(\vec\psi^m),~\mat{\bD}(\vec\chi)\bigr) + \bigl(\lambda^{m-1}\,\bG\cdot\vec\psi^m,~\bG\cdot\vec\chi\bigr)=0\qquad\forall\vec\chi\in\mathbb{Y}_0^{m-1},\label{eq:elastic}
\end{equation}
where we introduce 
\begin{equation}
\lambda^{m-1}_{|o_j^{m-1}} = 1 + \frac{\max\limits_{o\in\mathscr{T}^{m-1}} |o| - \min\limits_{o\in\mathscr{T}^{m-1}}|o|}{|o_j^{m-1}|},
\qquad j=1,\ldots,J_\mR,\nn
\end{equation}
to limit the distortion of small elements \cite{Masud1997space,Zhao2020energy}.

In view of \eqref{eq:meshQU}, it is natural to introduce the approximation of the mesh velocity $\vec w$ at time $t_m$ as 
\begin{equation}
\vec W^m(\vec x): = \sum_{k=1}^K\left(\frac{\vec q_k^{m} - \vec q_k^{m-1}}{\ttau}\right)\,\phi_k^m(\vec x),\quad \vec x\in\mR_\pm^m,
\label{eq:DMeshV}
\end{equation} 
where  $\phi_k^m(\cdot)$ is the nodal basis function of $S_1^m$ at $\vec q_k^m$. The corresponding discrete ALE mappings are given by 
\begin{align}
\mathcal{\vec A}^m[t](\vec x): = \vec\id - (t_m - t)\vec W^m(\vec x)=\sum_{k = 1}^{K}\left(\frac{t_{m} -t}{\ttau}\,\vec q_k^{m-1}+ \frac{t-t_{m-1}}{\ttau}\,\vec q_k^{m}\right)\phi_k^m(\vec x),\quad\forall t\in [t_{m-1},~t_{m}],\quad\vec x\in\mR,
\label{eq:DALE}
\end{align}
where $\vec\id$ is the identity function on $\mR$. It is easy to see that $\mathcal{\vec A}^m[t_m](\vec x)$ is the identity map, and 
\begin{equation}
\mathcal{A}^m[t_{m-1}](\vec x) = \vec\id - \ttau\,\vec W^m(\vec x)\quad\forall\vec x\in\mR_\pm^m
,\quad\mbox{with}\quad \mathcal{\vec A}^m[t_{m-1}](\mR_\pm^m) = \mR_\pm^{m-1}.
\end{equation}
The Jacobian determinant of the linear map $\mathcal{\vec A}^m[t_{m-1}]$ is given by
\begin{equation}
\mathcal{J}^{m}(\vec x): = {\rm det}(\bG\vec\id - \ttau\bG\vec W^m)= 1 - \ttau\,\bG\cdot\vec W^m + O(\ttau^2),\qquad \vec x\in\mR_\pm^m.\label{eq:Jacobian}
\end{equation}

We have the following lemma about the discrete ALE mappings. 
\begin{lem}\label{lem:ncALET}
Suppose that $\varphi\in L^2(\mR_\pm^{m})$. Then it holds that
\begin{equation}
\int_{\mR^m_\pm}\varphi\,\Bigl(r- \,\ttau\,[\vec W^m\cdot\vec e_1]\Bigr)\,\mathcal{J}^{m}\,\drz = \int_{\mR^{m-1}_\pm}\varphi\circ\mathcal{\vec A}^m[t_{m-1}]^{-1}\,r\,\drz.\label{eq:nALET}
\end{equation}
Moreover, let $\mR_\pm^h(t)=\mathcal{\vec A}^m[t](\mR_\pm^m)$ for $t\in[t_{m-1}, t_m]$. Then it holds 
\begin{equation}
\int_{\mR_\pm^m}\varphi\,r\,\drz - \int_{\mR_\pm^{m-1}}\varphi\circ\mathcal{\vec A}^m[t_{m-1}]^{-1}\,r\,\drz  = \int_{t_{m-1}}^{t_m}\int_{\mR_\pm^h(t)}\varphi\circ\mathcal{\vec A}^m[t]^{-1}\nabla\cdot[r\,\vec W^m\circ\mathcal{\vec A}^m[t]^{-1}]\,\drz\rd t.
\label{eq:cALET}
\end{equation}
\end{lem}
\begin{proof}
We recall the definition of $\mathcal{J}^m$  in \eqref{eq:Jacobian} as well as the fact that $\mathcal{\vec A}^m[t_{m-1}] = \vec\id  - \ttau\,\vec W^m$ and $\vec\id\cdot\vec e_1 = r$. Then using the change of variables $\vec x = \mathcal{\vec A}^m[t_{m-1}](\vec y)$, it is straightforward to obtain
\begin{equation}
\int_{\mR^m_\pm}\varphi\,\Bigl(\mathcal{\vec A}^m[t_{m-1}]\cdot\vec e_1\Bigr)\,\mathcal{J}^m\,\drz = \int_{\mR^{m-1}_\pm}\varphi\circ\mathcal{\vec A}^m[t_{m-1}]^{-1}\,r\,\drz, \nn
\end{equation}
which then implies \eqref{eq:nALET}.

It follows directly from \eqref{eq:DALE} that
\begin{equation}
\partial_t^\circ(\vec\chi\circ\mathcal{\vec A}^m[t]^{-1})=\vec 0\qquad\forall\vec\chi\in\bU^m.
\label{eq:zeroD}
\end{equation}
Then applying \eqref{eq:ddtALEframe} to the domain $\mR^h_\pm(t)$ and using \eqref{eq:zeroD} yields that
\begin{equation}
\ddt\int_{\mR_\pm^h(t)}\varphi\circ\mathcal{\vec A}^m[t]^{-1}\,r\,\drz  = \int_{\mR_\pm^h(t)}\varphi\circ\mathcal{\vec A}^m[t]^{-1}\nabla\cdot[r\,\vec W^m\circ\mathcal{\vec A}^m[t]^{-1}]\,\drz,\label{eq:cALE1}
\end{equation}
which leads to \eqref{eq:cALET} directly after integrating  with respect to $t$ from $t_{m-1}$ to $t_m$. 
\end{proof}

\begin{rem}
We note that \eqref{eq:cALET} can be regarded as the axisymmetric analogue of the geometric conservation law that was introduced in \cite[(3.22)]{Nobile1999}. Moreover, on applying the change of variables $\vec x = \mathcal{\vec A}^m[t]^{-1}(\vec y)$ to \eqref{eq:cALE1}, we obtain
\begin{align}
&\int_{\mR_\pm^h(t)}\varphi\circ\mathcal{\vec A}^m[t]^{-1}\nabla\cdot[r\,\vec W^m\circ\mathcal{\vec A}^m[t]^{-1}]\,\drz=\ddt\int_{\mR_\pm^h(t)}\varphi\circ\mathcal{\vec A}^m[t]^{-1}\,r\,\drz \nn\\
 &=\ddt\int_{\mR_\pm^m}\varphi\det\mat{G}(\vec x,t)\,(\mathcal{\vec A}^m[t]\cdot\vec e_1)\drz =\int_{\mR_\pm^m}\varphi\ddt\left\{\det\mat{G}(\vec x,t)\,(\mathcal{\vec A}^m[t]\cdot\vec e_1)\right\}\,\drz,
 \label{eq:cALE5}
\end{align}
where $\mat{G}(\vec x, t) = \frac{\partial\mathcal{\vec A}^m[t](\vec x)}{\partial\vec x}$. On recalling \eqref{eq:DALE}, we see that \eqref{eq:cALE5} is a polynomial of degree $2$ in the variable $t$. Thus the right hand side of \eqref{eq:cALET} can be integrated exactly with respect to $t$ via Simpson's rule.

\end{rem}

\subsection{Nonconservative ALE approximations}
\label{sec:nALEm}
In the following, we denote by $\vec U^m$, $P^m$, $\varkappa^m$ and $\kappa^m$ the numerical approximations of $\vec u(\cdot,t)$, $p(\cdot,t)$, $\varkappa(\cdot,t)$ and $\kappa(\cdot, t)$ at time $t_m$, respectively. We introduce 
\begin{equation}
\vec\chi_{_\mathcal{A}} = \vec\chi\circ\mathcal{\vec A}^m[t_{m-1}]\in\bU^m\qquad\mbox{for}\quad\vec \chi\in\bU^{m-1}.\label{eq:at}
\end{equation}

We first consider an unconditionally stable approximation of the weak formulation \eqref{eqn:weakform} as follows. Let $\Gamma^0:=\vec X^0(\cdot)\in V_\partial^h$ and $\vec U^0\in \bU^0$
be the approximations of the initial interface and velocity field, respectively. Moreover, we set $\Gamma^{-1}=\Gamma^0$, $\mR_\pm^{-1}=\mR_\pm^0$ with $\vec W^0=\vec 0$ and $\mathcal{J}^0(\vec x) = 1$. Then for $m\geq 0$, we seek $\vec U^{m+1}\in \bU^m$, 
$P^{m+1}\in\bP^m$, $\vec X^{m+1}\in V_\partial^h$ 
and $\varkappa^{m+1}\in V^h$ such that
\begin{subequations}\label{eqn:stabfd}
\begin{align}
\label{eq:stabfd1}
&\Bigl(\rho^m\frac{\vec U^{m+1} - \vec U^m_\mathcal{A}\,\sqrt{\Bigl(1 - \,\ttau\,r^{-1}[\vec W^m\cdot\vec e_1]\Bigr)\,\mathcal{J}^m}}{\ttau},~\vec\chi\,r\Bigr)^\diamond +\mathscr{A}(\rho^m,\vec U^m_\mathcal{A}-\vec W^m; \vec U^{m+1},\vec\chi)\nn\\
&\hspace{1.5cm}+2\Bigl(\mu^m\,r^{-1}[\vec U^{m+1}\cdot\vec e_1],~[\vec\chi\cdot\vec e_1]\Bigr)^\diamond+2\Bigl(\mu^m\,r\,\mat{\bD}(\vec U^{m+1}),~\mat{\bD}(\vec\chi)\Bigr)- \Bigl(P^{m+1},~\bG\cdot[r\vec\chi]\Bigr) \nn\\
&\hspace{1.5cm} -\,\gamma\,\Big\langle(\vec X^{m}\cdot\vec e_1)\,\varkappa^{m+1},~\vec\nu^m\cdot\vec\chi\,|\vec X^m_\alpha|\Big\rangle  = \Bigl(\rho^m\vec g,~\vec\chi\,r\Bigr)\qquad\forall\vec\chi\in\bU^m,\\[0.5em]
\label{eq:stabfd2}
&\hspace{0.1cm}\Bigl(\bG\cdot[r\vec U^{m+1}],~q\Bigr)=0\qquad\forall q\in \bP^m, \\[0.5em]
\label{eq:stabfd3}
&\frac{1}{\ttau}\Big\langle(\vec X^m\cdot\vec e_1)(\vec X^{m+1}-\vec X^m),~\zeta\,\vec\nu^m\,|\vec X_\alpha^m|\Big\rangle-\Big\langle(\vec X^m\cdot\vec e_1)\vec U^{m+1}\cdot\vec\nu^m,~\zeta\,|\vec X^m_\alpha|\Big\rangle=0\qquad\forall\zeta\in V^h,\\[0.5em]
&\Big\langle(\vec X^m\cdot\vec e_1)\varkappa^{m+1},~\vec\eta\cdot\vec\nu^m\,|\vec X_\alpha^m|\Big\rangle +\Big\langle\vec\eta\cdot\vec e_1,~|\vec X^{m+1}_\alpha|\Big\rangle + \Big\langle(\vec X^m\cdot\vec e_1)\vec X^{m+1}_\alpha,~\vec\eta_\alpha\,|\vec X^m_\alpha|^{-1}\Big\rangle=0\qquad\forall\vec\eta\in V_\partial^h,
\label{eq:stabfd4}
\end{align}
\end{subequations}
where $\vec U^m_{\mathcal{A}}$ and $\vec W^m$ are defined via \eqref{eq:at} and \eqref{eq:DMeshV}, respectively, and we set $\Gamma^{m+1}=\vec X^{m+1}(\bI)$ to construct the new bulk mesh through \eqref{eq:meshQU} and \eqref{eq:elastic}. Moreover, $(\cdot,\cdot)^\diamond$ represents an approximation of the inner product $(\cdot,\cdot)$ using a high-order Gauss quadrature rule that is exact for polynomial of degree at most $5$. 

We next show that \eqref{eq:stabfd1} is in fact a consistent discretization of \eqref{eq:weak1}. Recalling \eqref{eq:Jacobian} and using a Taylor expansion of the first term in \eqref{eq:stabfd1} yields 
\begin{align}
&\Bigl(\rho^m\frac{\vec U^{m+1} - \vec U^m_{\mathcal{A}}\,\sqrt{\Bigl(1 - \ttau\,r^{-1}[\vec W^m\cdot\vec e_1]\Bigr)\mathcal{J}^m}}{\ttau},~\vec\chi\,r\Bigr)^\diamond\nn\\
&\hspace{0.5cm} = \Bigl(\rho^m\frac{\vec U^{m+1}-\vec U^m_\mathcal{A}\,\sqrt{\Bigl(1 - \ttau\,r^{-1}[\vec W^m\cdot\vec e_1]\Bigr)\Bigl(1-\ttau\,\bG\cdot\vec W^m\Bigr) + O(\ttau^2)}}{\ttau},~\vec\chi\,r\Bigr)^\diamond\nn\\
&\hspace{0.5cm}=\Bigl(\rho^m\frac{\vec U^{m+1} - \vec U^m_\mathcal{A}}{\ttau},~\vec\chi\,r\Bigr)^\diamond + \frac{1}{2}\Bigl(\rho^m\bG\cdot\vec W^m,~\vec U^m_\mathcal{A}\cdot\vec\chi\,r\Bigr)^\diamond + \frac{1}{2}\Bigl(\rho^m[\vec W^m\cdot\vec e_1],~\vec U^m_\mathcal{A}\cdot\vec\chi\Bigr)^\diamond+ O(\ttau)\nn\\
&\hspace{0.5cm}=\Bigl(\rho^m\frac{\vec U^{m+1} - \vec U^m_\mathcal{A}}{\ttau},~\vec\chi\,r\Bigr)^\diamond + \frac{1}{2}\Bigl(\rho^m\bG\cdot[r\vec W^m],~\vec U^m_\mathcal{A}\cdot\vec\chi\Bigr)+ O(\ttau),\label{eq:conT}
\end{align}
which is a consistent temporal discretization of the first two terms in \eqref{eq:weak1}. The special approximation in \eqref{eq:conT} allows a stable discretization with a decreasing discrete kinetic energy. We also note that \eqref{eqn:stabfd} is nonlinear due to the presence of $|\vec X^{m+1}_\alpha|$ in \eqref{eq:stabfd4}, which contributes to the stability of the interface energy. As a consequence, we have the following theorem for the introduced method \eqref{eqn:stabfd}, which mimics the energy stability \eqref{eq:energylaw} on the discrete level.
\begin{thm}\label{thm:ES} Let $(\vec U^{m+1}, P^{m+1}, \vec X^{m+1}, \varkappa^{m+1})$ be a solution to \eqref{eqn:stabfd} for $m=0,1,\cdots,M-1$. Then it holds that 
\begin{align}
&\frac{1}{2\pi}E(\rho^k,\vec U^{k+1},\vec X^{k+1}) +  2\ttau\sum_{m=0}^k\Bigl(\norm{\sqrt{\mu^m\,r^{-1}}\,(\vec U^{m+1}\cdot\vec e_1)}_\diamond^2 + \norm{\sqrt{\mu^m\,r}\,\mat{\bD}(\vec U^{m+1})}^2\Bigr)\nn\\
&\hspace{1.5cm}\leq \frac{1}{2\pi}E(\rho^{0},\vec U^0, \vec X^0)+ \ttau\sum_{m=0}^k\bigl(\rho^m\,r\,\vec g,~\vec U^{m+1}\bigr),\qquad k=0,1,\cdots, M-1,
\label{eq:DES}
\end{align}
where $\norm{\cdot}$ and $\norm{\cdot}_\diamond$ are the induced norms of the inner products $(\cdot,\cdot)$ and $(\cdot,\cdot)^\diamond$, respectively.
\end{thm}
\begin{proof}
Setting $\vec\chi= \ttau\vec U^{m+1}$ in \eqref{eq:stabfd1}, $q = P^{m+1}$ in \eqref{eq:stabfd2}, $\zeta = \ttau\gamma\varkappa^{m+1}$ in \eqref{eq:stabfd3} and $\vec\eta = \gamma(\vec X^{m+1} - \vec X^m)$ in \eqref{eq:stabfd4}, and combining these four equations, yields 
\begin{align}
&\Bigl(\rho^m\, \delta\vec U^m,~\vec U^{m+1}\,r\Bigr) + 2\ttau\,\Bigl(\mu^m\,r^{-1}[\vec U^{m+1}\cdot\vec e_1],~[\vec U^{m+1}\cdot\vec e_1]\Bigr) + 2\ttau\Bigl(\mu^m\,r\,\mat{\bD}(\vec U^{m+1}),~\mat{\bD}(\vec U^{m+1})\Bigr)\nn\\
&\hspace{1cm}+\gamma\Big\langle(\vec X^{m+1}-\vec X^m)\cdot\vec e_1,~|\vec X_\alpha^{m+1}|\Big\rangle + \gamma\Big\langle(\vec X^m\cdot\vec e_1)\vec X^{m+1}_\alpha,~(\vec X^{m+1}-\vec X^m)_\alpha\,|\vec X^m_\alpha|^{-1}\Big\rangle  = \ttau\Big(\rho^m\,r\,\vec g,~\vec U^{m+1}\Bigr),\label{eq:AEstab}
\end{align}
where we denote $\delta\vec U^m = \vec U^{m+1}-\vec U_\mathcal{A}^m\sqrt{\Bigl(1 - \ttau\,r^{-1}[\vec W^m\cdot\vec e_1]\Bigr)\mathcal{J}^m}$.

By the inequality $\vec a\cdot(\vec a- \vec b)\geq |\vec b|(|\vec a| - |\vec b|)$, we have
\begin{align}
&\Big\langle(\vec X^{m+1}-\vec X^m)\cdot\vec e_1,~|\vec X_\alpha^{m+1}|\Big\rangle + \Big\langle(\vec X^m\cdot\vec e_1)\vec X^{m+1}_\alpha,~(\vec X^{m+1}-\vec X^m)_\alpha\,|\vec X^m_\alpha|^{-1}\Big\rangle \nn\\
&\hspace{1cm}\geq \Big\langle(\vec X^{m+1}-\vec X^m)\cdot\vec e_1,~|\vec X_\alpha^{m+1}|\Big\rangle + \Big\langle(\vec X^m\cdot\vec e_1),~|\vec X^{m+1}_\alpha|-|\vec X^m_\alpha|)\Big\rangle\nn\\
&\hspace{1cm}= \Big\langle(\vec X^{m+1}\cdot\vec e_1),~|\vec X_\alpha^{m+1}|\Big\rangle - \Big\langle(\vec X^m\cdot\vec e_1),~|\vec X^m_\alpha|\Big\rangle=\frac{1}{2\pi}\left\{A(\vec X^{m+1}) - A(\vec X^m)\right\},\label{eq:Astab}
\end{align}
where we invoke the definition of $A(\vec X^m)$ in \eqref{eq:area}.  

Let $(\cdot,\cdot)_{\mR_\pm^m}$ denote the $L^2$--inner product over $\mR_{\pm}^m$. Moreover, we let $(\cdot,\cdot)_{\mR_\pm^m}^\diamond$ be an approximation of $(\cdot,\cdot)_{\mR_\pm^m}$ using the prescribed high-order Gauss quadrature rule. Using the inequality $2\vec a\cdot(\vec a -\vec b)\geq |\vec a|^2 - |\vec b|^2$ yields that
\begin{align}\label{eq:kEstab1}
\Bigl(\delta\vec U^m,~\vec U^{m+1}\,r\Bigr)_{\mR_\pm^m}^\diamond  &\geq \frac{1}{2}\Bigl(|\vec U^{m+1}|^2,~r\Bigr)_{\mR_\pm^m}^\diamond - \frac{1}{2}\Bigl(|\vec U^m_\mathcal{A}|^2\,\bigl(r -\ttau\,[\vec W^m\cdot\vec e_1]\bigr)\,\mathcal{J}^m, ~1\Bigr)_{\mR_\pm^m}^\diamond\nn\\
&= \frac{1}{2}\Bigl(|\vec U^{m+1}|^2,~r\Bigr)_{\mR_\pm^m}  - \frac{1}{2}\Bigl(|\vec U^m|^2, ~r\Bigr)_{\mR_\pm^{m-1}},
\end{align}
where for the last equality we applied \eqref{eq:nALET} with $\varphi = |\vec U_{\mathcal{A}}^m|^2$. We then multiply \eqref{eq:kEstab1} with $\rho_\pm$ and combine the two equations to give
\begin{equation}\label{eq:kEstab}
\Bigl(\rho^m\,\delta\vec U^m,~\vec U^{m+1}\,r\Bigr)^\diamond\geq \frac{1}{2}\Bigl(\rho^m\,|\vec U^{m+1}|^2,~r\Bigr) - \frac{1}{2}\Bigl(\rho^{m-1}\,|\vec U^m|^2,~r\Bigr),
\end{equation}
on recalling \eqref{eq:Ddensity}. Inserting \eqref{eq:Astab} and \eqref{eq:kEstab} into \eqref{eq:AEstab}, and recalling \eqref{eq:Energy} leads to 
\begin{align}
&\frac{1}{2\pi}E(\rho^m, \vec U^{m+1}, \vec X^{m+1}) + 2\ttau\Bigl(\norm{\sqrt{\mu^m\,r^{-1}}\,(\vec U^{m+1}\cdot\vec e_1)}_\diamond^2 + \norm{\sqrt{\mu^m\,r}\,\mat{\bD}(\vec U^{m+1})}^2\Bigr)\nn\\
&\hspace{1cm}\leq \frac{1}{2\pi}E(\rho^{m-1},\vec U^m, \vec X^m)+ \ttau\bigl(\rho^m\,r\,\vec g,~\vec U^{m+1}\bigr).\label{eq:Astabstep}
\end{align}
Summing \eqref{eq:Astabstep} for $m=0,\cdots,k$ yields \eqref{eq:DES} immediately on recalling $\mR_\pm^{-1}=\mR_\pm^0$.
\end{proof}

We next consider a linear approximation of \eqref{eqn:Aweakform} as follows. With the discrete initial data as before, for $m\geq 0$, we find $\vec U^{m+1}\in \bU^m$, 
$P^{m+1}\in\bP^m$, $\vec X^{m+1}\in V_\partial^h$ 
and $\kappa^{m+1}\in V^h$ such that
\begin{subequations}\label{eqn:meshfd}
\begin{align}
\label{eq:meshfd1}
&\Bigl(\rho^m\frac{\vec U^{m+1} - \vec U^m_\mathcal{A}\,\sqrt{\Bigl(1 - \,\ttau\,r^{-1}[\vec W^m\cdot\vec e_1]\Bigr)\,\mathcal{J}^m}}{\ttau},~\vec\chi\,r\Bigr)^\diamond +\mathscr{A}(\rho^m,\vec U^m_\mathcal{A}-\vec W^m; \vec U^{m+1},\vec\chi)\nn\\
&\hspace{1.5cm}+2\Bigl(\mu^m\,r^{-1}[\vec U^{m+1}\cdot\vec e_1],~[\vec\chi\cdot\vec e_1]\Bigr)^\diamond+2\Bigl(\mu^m\,r\,\mat{\bD}(\vec U^{m+1}),~\mat{\bD}(\vec\chi)\Bigr)- \Bigl(P^{m+1},~\bG\cdot[r\vec\chi]\Bigr) \nn\\
&\hspace{1.5cm} -\,\gamma\,\Big\langle(\vec X^{m}\cdot\vec e_1)\,\kappa^{m+1}-\vec\nu^m\cdot\vec e_1,~\vec\nu^m\cdot\vec\chi\,|\vec X^m_\alpha|\Big\rangle  = \Bigl(\rho^m\vec g,~\vec\chi\,r\Bigr)\qquad\forall\vec\chi\in\bU^m,\\[0.5em]
\label{eq:meshfd2}
&\hspace{0.1cm}\Bigl(\bG\cdot[r\vec U^{m+1}],~q\Bigr)=0\qquad\forall q\in \bP^m, \\[0.5em]
\label{eq:meshfd3}
&\frac{1}{\ttau}\Big\langle(\vec X^m\cdot\vec e_1)(\vec X^{m+1}-\vec X^m),~\zeta\,\vec\nu^m\,|\vec X_\alpha^m|\Big\rangle-\Big\langle(\vec X^m\cdot\vec e_1)\vec U^{m+1}\cdot\vec\nu^m,~\zeta\,|\vec X^m_\alpha|\Big\rangle=0\qquad\forall\zeta\in V^h,\\[0.5em]
&\Big\langle\kappa^{m+1}\vec\nu^m,~\vec\eta\,|\vec X_\alpha^m|\Big\rangle^h  + \Big\langle\vec X^{m+1}_\alpha,~\vec\eta_\alpha\,|\vec X^m_\alpha|^{-1}\Big\rangle=0\qquad\forall\vec\eta\in V_\partial^h,
\label{eq:meshfd4}
\end{align}
\end{subequations}
where we introduced $\langle\cdot,\cdot\rangle^h$ as the mass-lumped $L^2$--inner product over $\bI$, i.e., 
\[
\big\langle \vec v, \vec u \big\rangle^h = \frac{1}{2}\,h\sum_{j=1}^{J_\Gamma} 
\left[(\vec v\cdot\vec u)(\alpha_j^-) + (\vec v \cdot \vec u)(\alpha_{j-1}^+)\right]\quad\mbox{with}\quad g(\alpha_j^\pm)=\underset{\delta\searrow 0}{\lim}\ g(\alpha_j\pm\delta),
\]
for two piecewise continuous functions $\vec v, \vec w$. 

It does not seem possible to obtain discrete energy stability for the above introduced method \eqref{eqn:meshfd}. Nevertheless, \eqref{eqn:meshfd} gives rise to a system of linear equations, and the approximation in \eqref{eq:meshfd4} leads to equidistribution, meaning that the mesh points on $\Gamma^m$ tend to be distributed at evenly spaced arc length. The interested reader is referred to \cite{BGN07, Barrett20} for detailed discussions of this property.

\subsection{Conservative ALE approximations}\label{sec:cALEm}

Based on \eqref{eq:weakc}, we propose a conservative ALE method which satisfies an unconditional stability estimate as well. With the same discrete initial data as before, for $m\geq 0$, we find $\vec U^{m+1}\in \bU^m$, 
$P^{m+1}\in\bP^m$, $\vec X^{m+1}\in V_\partial^h$
and $\varkappa^{m+1}\in V^h$ such that
\begin{subequations}
\label{eqn:cstabfd}
\begin{align}
\label{eq:cstabfd1}
&\frac{1}{\ttau}\left[\Bigl(\rho^m\vec U^{m+1},~\vec\chi\,r\Bigr) - \Bigl(\rho^{m-1}\vec U^m,~\vec\chi\circ\mathcal{\vec A}^m[t_{m-1}]^{-1}\,r\Bigr)\right]- \mathscr{B}(\rho^m,\vec W^m;\vec U^{m+1}\cdot\vec\chi)+\mathscr{A}(\rho^m,\vec U^m_\mathcal{A}-\vec W^m; \vec U^{m+1},\vec\chi)\nn\\
&\hspace{1.5cm}+2\Bigl(\mu^m\,r^{-1}[\vec U^{m+1}\cdot\vec e_1],~[\vec\chi\cdot\vec e_1]\Bigr)^\diamond+2\Bigl(\mu^m\,r\,\mat{\bD}(\vec U^{m+1}),~\mat{\bD}(\vec\chi)\Bigr)- \Bigl(P^{m+1},~\bG\cdot[r\vec\chi]\Bigr) \nn\\
&\hspace{1.5cm} -\,\gamma\,\Big\langle(\vec X^{m}\cdot\vec e_1)\,\varkappa^{m+1},~\vec\nu^m\cdot\vec\chi\,|\vec X^m_\alpha|\Big\rangle  = \Bigl(\rho^m\vec g,~\vec\chi\,r\Bigr)\qquad\forall\vec\chi\in\bU^m,\\
\label{eq:cstabfd2}
&\hspace{0.1cm}\Bigl(\bG\cdot[r\vec U^{m+1}],~q\Bigr)=0\qquad\forall q\in \bP^m, \\[0.5em]
\label{eq:cstabfd3}
&\frac{1}{\ttau}\Big\langle(\vec X^m\cdot\vec e_1)(\vec X^{m+1}-\vec X^m),~\zeta\,\vec\nu^m\,|\vec X_\alpha^m|\Big\rangle-\Big\langle(\vec X^m\cdot\vec e_1)\vec U^{m+1}\cdot\vec\nu^m,~\zeta\,|\vec X^m_\alpha|\Big\rangle=0\qquad\forall\zeta\in V^h,\\[0.5em]
&\Big\langle(\vec X^m\cdot\vec e_1)\varkappa^{m+1},~\vec\eta\cdot\vec\nu^m\,|\vec X_\alpha^m|\Big\rangle +\Big\langle\vec\eta\cdot\vec e_1,~|\vec X^{m+1}_\alpha|\Big\rangle + \Big\langle(\vec X^m\cdot\vec e_1)\vec X^{m+1}_\alpha,~\vec\eta_\alpha\,|\vec X^m_\alpha|^{-1}\Big\rangle=0\qquad\forall\vec\eta\in V_\partial^h,
\label{eq:cstabfd4}
\end{align}
\end{subequations}
where we introduced the time-integrated term 
\begin{equation}
\mathscr{B}(\rho^m,\vec W^m;\varphi) = \frac{1}{2\ttau}\int_{t_m}^{t_{m+1}}\Bigl(\rho^m\circ\mathcal{\vec A}^m[t]^{-1}\nabla\cdot[r\,\vec W^m\circ\mathcal{\vec A}^m[t]^{-1}],~\varphi\circ\mathcal{\vec A}^m[t]^{-1}\Bigr)\,\rd t.\nn
\end{equation}

We have the following theorem for the introduced method \eqref{eqn:cstabfd}.
\begin{thm}\label{thm:cES}
If $(\vec U^{m+1}, P^{m+1}, \vec X^{m+1}, \varkappa^{m+1})$ is a solution to \eqref{eqn:cstabfd} for $m=0,\cdots,M-1$,  then the energy stability estimate \eqref{eq:DES} holds.
\end{thm}
\begin{proof}
We choose $\vec\chi= \ttau\vec U^{m+1}$ in \eqref{eq:cstabfd1}, $q = P^{m+1}$ in \eqref{eq:cstabfd2}, $\zeta = \ttau\gamma\varkappa^{m+1}$ in \eqref{eq:cstabfd3} and $\vec\eta = \gamma(\vec X^{m+1} - \vec X^m)$ in \eqref{eq:cstabfd4} and combine these equations.  Then the proof is very similar to that of Theorem \ref{thm:ES}, and it remains to show the stability for the fluid kinetic energy 
\begin{align}
&\Bigl(\rho^m\vec U^{m+1},~\vec U^{m+1}\,r\Bigr)-\Bigl(\rho^{m-1}\vec U^m,~\vec U^{m+1}\circ\mathcal{A}^m[t_{m-1}]^{-1}\,r\Bigr)-\ttau\,\mathscr{B}(\rho^m,\vec W^m;|\vec U^{m+1}|^2)\nn\\
&\hspace{1cm}\geq \frac{1}{2}\Bigl(\rho^m\,|\vec U^{m+1}|^2, r\Bigr) - \frac{1}{2}\Bigl(\rho^{m-1}|\vec U^m|^2,~r\Bigr).\label{eq:cALE2}
\end{align}
 
In view of the inequality $2\vec a\cdot\vec b\leq |\vec a|^2 +|\vec b|^2$, we have
\begin{equation}
\frac{1}{2}\Bigl(\rho^{m-1}|\vec U^m|^2,~r\Bigr) + \frac{1}{2}\Bigl(\rho^{m-1}|\vec U^{m+1}\circ\mathcal{\vec A}^m[t_{m-1}]^{-1}|^2,~r\Bigr)\geq \Bigl(\rho^{m-1}\vec U^m,~\vec U^{m+1}\circ\mathcal{\vec A}^m[t_{m-1}]^{-1}\,r\Bigr).\label{eq:cALE3}
\end{equation}
Moreover, choosing $\varphi=|\vec U^{m+1}|^2$ in \eqref{eq:cALET}, then multiplying \eqref{eq:cALET} with $\rho_\pm$ and adding the two equations leads to
\begin{equation}
\frac{1}{2}\Bigl(\rho^m\,|\vec U^{m+1}|^2,~r\Bigr) -\frac{1}{2} \Bigl(\rho^{m-1}\,|\vec U^{m+1}\circ\mathcal{\vec A}^m[t_{m-1}]^{-1}|^2,~r\Bigr)=\ttau\,\mathscr{B}(\rho^m,\vec W^m; |\vec U^{m+1}|^2).\label{eq:cALE4}
\end{equation}
Combining \eqref{eq:cALE3} and \eqref{eq:cALE4} then implies \eqref{eq:cALE2}. 
\end{proof}

\medskip
Based on \eqref{eq:Aweakc} and \eqref{eqn:meshfd}, it is natural to introduce the following conservative ALE method that enjoys the property of equidistribution. With the same discrete initial data, for $m\geq 0$, we find $\vec U^{m+1}\in \bU^m$, 
$P^{m+1}\in\bP^m$, $\vec X^{m+1}\in V_\partial^h$ 
and $\kappa^{m+1}\in V^h$ such that
\begin{subequations}\label{eqn:cmeshfd}
\begin{align}
\label{eq:cmeshfd}
&\frac{1}{\ttau}\left[\Bigl(\rho^m\vec U^{m+1},~\vec\chi\,r\Bigr) - \Bigl(\rho^{m-1}\vec U^m,~\vec\chi\circ\mathcal{\vec A}^m[t_{m-1}]^{-1}\,r\Bigr)\right]- \mathscr{B}(\rho^m,\vec W^m;\vec U^{m+1}\cdot\vec\chi)+\mathscr{A}(\rho^m,\vec U^m_\mathcal{A}-\vec W^m; \vec U^{m+1},\vec\chi)\nn\\
&\hspace{1.5cm}+2\Bigl(\mu^m\,r^{-1}[\vec U^{m+1}\cdot\vec e_1],~[\vec\chi\cdot\vec e_1]\Bigr)^\diamond+2\Bigl(\mu^m\,r\,\mat{\bD}(\vec U^{m+1}),~\mat{\bD}(\vec\chi)\Bigr)- \Bigl(P^{m+1},~\bG\cdot[r\vec\chi]\Bigr) \nn\\
&\hspace{1.5cm} -\,\gamma\,\Big\langle(\vec X^{m}\cdot\vec e_1)\,\kappa^{m+1} - \vec\nu^m\cdot\vec e_1,~\vec\nu^m\cdot\vec\chi\,|\vec X^m_\alpha|\Big\rangle  = \Bigl(\rho^m\vec g,~\vec\chi\,r\Bigr)\qquad\forall\vec\chi\in\bU^m,\\
\label{eq:cmeshfd2}
&\hspace{0.1cm}\Bigl(\bG\cdot[r\vec U^{m+1}],~q\Bigr)=0\qquad\forall q\in \bP^m, \\[0.5em]
\label{eq:cmeshfd3}
&\frac{1}{\ttau}\Big\langle(\vec X^m\cdot\vec e_1)(\vec X^{m+1}-\vec X^m),~\zeta\,\vec\nu^m\,|\vec X_\alpha^m|\Big\rangle-\Big\langle(\vec X^m\cdot\vec e_1)\vec U^{m+1}\cdot\vec\nu^m,~\zeta\,|\vec X^m_\alpha|\Big\rangle=0\qquad\forall\zeta\in V^h,\\[0.5em]
&\Big\langle\kappa^{m+1}\vec\nu^m,~\vec\eta\,|\vec X_\alpha^m|\Big\rangle^h  + \Big\langle\vec X^{m+1}_\alpha,~\vec\eta_\alpha\,|\vec X^m_\alpha|^{-1}\Big\rangle=0\qquad\forall\vec\eta\in V_\partial^h.
\label{eq:cmeshfd4}
\end{align}
\end{subequations}

\subsection{Volume-preserving approximations}\label{sec:vp}

We recall from \S\ref{sec:ALEWF} that it is possible to prove volume
preservation for the weak formulations on choosing suitable test functions, see \eqref{eq:vollaw}. However, for the ALE methods introduced in \S\ref{sec:nALEm} and \S\ref{sec:cALEm},  the volume of the two discrete phases will in general not be exactly preserved.
In fact, in order to enable an exact volume preservation on the discrete level, it turns out that suitable time-weighted approximations of the interface normals are necessary \cite{BZ21SPFEM,BGNZ22volume, GNZ23asy}. We have a discrete analogue of \eqref{eq:ddtM} as follows, and its proof can be found in \cite[Lemma 3.1]{BGNZ22volume}.
\begin{lem}\label{lem:vc} Let $\vec X^m\in V_\partial^h$ and $\vec X^{m+1}\in V_\partial^h$, then it holds that for $m=0,1,\cdots, M-1$
\begin{equation}
\vol(\vec X^{m+1})-\vol(\vec X^m) = 2\pi\big\langle\vec X^{m+1}-\vec X^m,~\vec f^{m+\frac{1}{2}}\big\rangle,
\end{equation}
where we introduced $\vec f^{m+\frac{1}{2}}\in [L^\infty(\bI)]^2$ as an appropriate treatment of the quantity $\vec f =(\vec\mZ\cdot\vec e_1)\,|\vec\mZ_\alpha|\,\vec\nu$:
\begin{align}
\label{eq:averagenor}
\vec f^{m+\frac{1}{2}}= -\frac{1}{6}\bigl[2(\vec X^m\cdot\vec e_1)\,\vec X^m_{\alpha}+2(\vec X^{m+1}\cdot\vec e_1)\,\vec X^{m+1}_{\alpha} + (\vec X^m\cdot\vec e_1)\,\vec X^{m+1}_{\alpha} + (\vec X^{m+1}\cdot\vec e_1)\,\vec X_{\alpha}^m\bigr]^\perp.
\end{align}
\end{lem}

Similarly to \cite{GNZ23asy}, we are now ready to adapt the methods \eqref{eqn:stabfd} and \eqref{eqn:cstabfd} to achieve structure-preserving approximations, meaning that the volume preservation and energy stability are satisfied on the discrete level. This can be done easily by replacing \eqref{eq:stabfd3}-\eqref{eq:stabfd4} and \eqref{eq:cstabfd3}--\eqref{eq:cstabfd4} with 
\begin{subequations}\label{eqn:spfd}
\begin{align}
\label{eq:spfd3}
&\frac{1}{\ttau}\big<\vec X^{m+1}-\vec X^m,~\zeta\,\vec f^{m+\frac{1}{2}}\big>-\big<(\vec X^m\cdot\vec e_1)\vec U^{m+1}\cdot\vec\nu^m,~\zeta\,|\vec X^m_\alpha|\big>=0\qquad\forall\zeta\in V^h,\\[0.5em]
&\big<\varkappa^{m+1},~\vec\eta\cdot\vec f^{m+\frac{1}{2}}\big> +\big<\vec\eta\cdot\vec e_1,~|\vec X^{m+1}_\alpha|\big> + \big<(\vec X^m\cdot\vec e_1)\vec X^{m+1}_\alpha,~\vec\eta_\alpha\,|\vec X^m_\alpha|^{-1}\big>=0\qquad\forall\vec\eta\in V_\partial^h.
\label{eq:spfd4}
\end{align}
\end{subequations}
We have the following theorem for the new adapted methods. 
\begin{thm}\label{thm:ESVC} Let $(\vec U^{m+1}, P^{m+1}, \vec X^{m+1}, \varkappa^{m+1})$ be a solution to the adapted methods of \eqref{eqn:stabfd} or \eqref{eqn:cstabfd} by replacing \eqref{eq:stabfd3}--\eqref{eq:stabfd4} or \eqref{eq:cstabfd3}--\eqref{eq:cstabfd4} with \eqref{eq:spfd3}--\eqref{eq:spfd4}. Then the energy stability estimate \eqref{eq:DES} holds.
Moreover, it holds for $m=0,\cdots, M-1$ that
\begin{equation}
\vol(\vec X^{m+1}) = \vol(\vec X^m).
\label{eq:DVC}
\end{equation}
\end{thm}
\begin{proof}
The stability estimate can be established in a similar manner to Theorem \ref{thm:ES} and Theorem \ref{thm:cES}. For the volume preservation, on recalling \eqref{eq:P2P1P0}, we choose 
\begin{equation}
q = (\mX_{_{\mR_-^m}}-\omega^m)\in\bP^m\quad\mbox{with}\qquad \omega^m = \frac{\int_{\mR_-^m}r\drz}{\int_\mR r\drz}\nn
\end{equation}
in \eqref{eq:stabfd2} or \eqref{eq:cstabfd2} and obtain 
\begin{align}
0&=\big(\bG\cdot[r\vec U^{m+1}],~q\big)=\int_{\mR_-^m}\bG\cdot[r\vec U^{m+1}]\,\drz - \omega^m\int_\mR\bG\cdot[r\vec U^{m+1}]\,\drz \nn\\
&=\int_{\Gamma^m}(\vec X^m\cdot\vec e_1)(\vec U^{m+1}\cdot\vec\nu^m)\drz=\Big\langle\vec X^m\cdot\vec e_1,~\vec U^{m+1}\cdot\vec\nu^m\,|\vec X^m_\alpha|\Big\rangle.\label{eq:ALEvc1}
\end{align}
On the other hand, setting $\zeta =\ttau$ in \eqref{eq:spfd3}  yields that
\begin{equation}
\Big\langle\vec X^{m+1}-\vec X^m,~\vec f^{m+\frac{1}{2}}\Big\rangle =\ttau\Big\langle\vec X^m\cdot\vec e_1,~\vec U^{m+1}\cdot\vec\nu^m\,|\vec X^m_\alpha|\Big\rangle= 0,\nn
\end{equation}
which implies \eqref{eq:DVC} on recalling \eqref{eq:ALEvc1}  and Lemma \ref{lem:vc}.
\end{proof}

Similarly, we have the following theorem for the volume-preserving variants of the methods \eqref{eqn:meshfd} and \eqref{eqn:cmeshfd}.
\begin{thm}\label{thm:meshVC}
Let $(\vec U^{m+1}, P^{m+1}, \vec X^{m+1}, \kappa^{m+1})$ be a solution to the adapted method of \eqref{eqn:meshfd} or \eqref{eqn:cmeshfd}, which is obtained by replacing \eqref{eq:meshfd3} or \eqref{eq:cmeshfd3} with \eqref{eq:spfd3}. Then it holds for $m=0,\cdots, M-1$ that
\begin{equation}
\vol(\vec X^{m+1}) = \vol(\vec X^m).\label{eq:DVCC}
\end{equation}
\end{thm}

\section{Numerical results}\label{sec:NR}

\setcounter{equation}{0}

\subsection{Solutions of the discrete systems}

The introduced ALE methods in \S\ref{sec:FEM} are summarized in Table \ref{tab:ALEsum}, where the adapted variants refer to the corresponding volume-preserving methods, indicated in the name of the scheme with a ``V''. We construct the discrete ALE mappings explicitly as described in \ref{sec:dale}, and the methods n-$\Equi^h$ and c-$\Equi^h$ lead to a system of linear equations, which can be written in matrix form as
\begin{equation}\label{eq:linsys}
\begin{pmatrix}
\mat{B}_{\mR}^m &\vec C_{_\mR}^m &-\gamma\vec N_{_\Gamma,_\mR}^m &0 \\[0.2em]
[\vec C_{_\mR}^m]^T & 0 &0 &0\\[0.2em]
[\vec N_{_\Gamma,_\mR}^m]^T &0 &0 &-\frac{1}{\ttau}[\vec N_{_\Gamma}^m]^T\\[0.2em]
0 &0 &\mathcal{\vec N}_{_\Gamma}^m &\mat{A}_{\Gamma}^m
\end{pmatrix}
\begin{pmatrix}
\vec U^{m+1}\\[0.2em]
P^{m+1} \\[0.2em]
\kappa^{m+1} \\[0.2em]
\delta\vec X^{m+1}
\end{pmatrix}=
\begin{pmatrix}
\vec c^m\\[0.2em]
0 \\[0.2em]
0 \\[0.2em]
-\mat{A}_\Gamma^m\vec X^{m}
\end{pmatrix}\,,
\end{equation}
where with a slight abuse of notations we denote by $\bigl(\vec U^{m+1},
P^{m+1},
\kappa^{m+1},
\delta\vec X^{m+1}\bigr)$ 
the coefficients of these finite element functions with respect to the standard bases of the corresponding finite element spaces, see \cite[(5.1a)]{BGN15stable}, and $\delta\vec X^{m+1}=\vec X^{m+1}-\vec X^m$. Denote 
\begin{equation}
\Xi_\Gamma^m:= \begin{pmatrix}
 0 & - \frac1{\ttau}\,\vec [N_\Gamma^m]^T \\[0.2em]
\mathcal{\vec N}_\Gamma^m & \mat{A}_\Gamma^m
\end{pmatrix} \,.\nn
\end{equation}
Applying the Schur complement approach to \eqref{eq:linsys} gives rise to two linear subsystems \cite{BGN15stable, Agnese20}
\begin{equation}\label{eq:subsys1}
\begin{pmatrix}
\mat{B}_{_\mR}^m + \gamma\begin{pmatrix}\vec N_{_\Gamma,_\mR}^m &0 \end{pmatrix} [\Xi_\Gamma^m]^{-1}\begin{pmatrix}[\vec N_{_\Gamma,_\mR}^m]^T\\[0.2em]
0\end{pmatrix} &\vec C_{_\mR}^m\\[0.4em]
[\vec C_{_\mR}^m]^T & 0
\end{pmatrix}
\begin{pmatrix}
\vec U^{m+1}\\[0.8em]
P^{m+1}
\end{pmatrix}
=\begin{pmatrix}
\vec c^m -  \gamma\begin{pmatrix}\vec N_{_\Gamma,_\mR}^m &0 \end{pmatrix} [\Xi_\Gamma^m]^{-1}\begin{pmatrix}
0\\[0.2em]
\mat{A}_{_\Gamma}^m\vec X^m
\end{pmatrix}
\\[0.4em]
0
\end{pmatrix}\,,
\end{equation}
and 
\begin{equation}\label{eq:subsys2}
\begin{pmatrix}
\kappa^{m+1}\\[0.2em]
\delta\vec X^{m+1}
\end{pmatrix}
=[\Xi_\Gamma^m]^{-1}
\begin{pmatrix}
-\vec N_{_\Gamma,_\mR}^m\vec U^{m+1}\\[0.2em]
-\mat{A}_{_\Gamma}^m\vec X^m
\end{pmatrix}\,.
\end{equation}
In practice, we solve \eqref{eq:subsys1} using the preconditioner GMRES method, and \eqref{eq:subsys2} with the help of a sparse LU factorization. 

\begin{table}[!htp]
\centering
\def\temptablewidth{0.75\textwidth}
\vspace{-0.5pt}
\caption{The properties of the introduced ALE methods and their volume-preserving variants.}
\vspace{0.2em}
\renewcommand\arraystretch{1.25}
\begin{tabular}{c|cccc}
\toprule[1pt]
 &\eqref{eqn:stabfd} &\eqref{eqn:meshfd} &\eqref{eqn:cstabfd} &\eqref{eqn:cmeshfd}\\ \hline
original method  &n-$\Stab^h$ & n-$\Equi^h$ & c-$\Stab^h$ &c-$\Equi^h$ \\ \hline 
unconditional stability &$\checkmark$ &$\times$ &$\checkmark$ &$\times$ \\ \hline
equidistribution &$\times$ &$\checkmark$ &$\times$ &$\checkmark$ \\ \hline
linear or nonlinear &nonlinear &linear &nonlinear &linear\\ \hline 
volume preservation &$\times$ &$\times$ &$\times$ &$\times$ \\ \midrule[1pt]\midrule[1pt]
adapted variants  &n-$\StabV^h$ & n-$\EquiV^h$ & c-$\StabV^h$ &c-$\EquiV^h$ \\ \hline 
unconditional stability &$\checkmark$ &$\times$ &$\checkmark$ &$\times$ \\ \hline
equidistribution &$\times$ &$\checkmark$ &$\times$ &$\checkmark$ \\ \hline
linear or nonlinear &nonlinear &polynomial &nonlinear &polynomial\\ \hline 
volume preservation &$\checkmark$ &$\checkmark$ &$\checkmark$ &$\checkmark$\\ 
\bottomrule[1pt]
\end{tabular}
\label{tab:ALEsum}
\end{table}

The nonlinear systems of equations resulting from the other introduced methods are solved via a lagged Picard-type iteration. For example, for the numerical method n-$\StabV^h$, which is given by \eqref{eq:stabfd1}--\eqref{eq:stabfd2} and \eqref{eq:spfd3}--\eqref{eq:spfd4}, the iterations at the $m$-th time step are given as follows. On setting $\vec X^{m+1,0} = \vec X^m$, for $\ell\geq 0$, we find $\vec U^{m+1,\ell+1}\in\bU^m$, $P^{m+1,\ell+1}\in\bP^m$, $\varkappa^{m+1,\ell+1}\in V^h$ and $\vec X^{m+1,\ell+1}\in V^h_\partial$ such that
\begin{subequations}\label{eqn:picard}
\begin{align}
\label{eq:picard1}
&\Bigl(\rho^m\frac{\vec U^{m+1,\ell+1} - \vec U^m_\mathcal{A}\,\sqrt{\Bigl(1 - \,\ttau\,r^{-1}[\vec W^m\cdot\vec e_1]\Bigr)\,\mathcal{J}^m}}{\ttau},~\vec\chi\,r\Bigr)^\diamond +\mathscr{A}(\rho^m,\vec U^m_\mathcal{A}-\vec W^m; \vec U^{m+1,\ell+1},\vec\chi)\nn\\
&\hspace{1.5cm}+2\Bigl(\mu^m\,r^{-1}[\vec U^{m+1,\ell+1}\cdot\vec e_1],~[\vec\chi\cdot\vec e_1]\Bigr)^\diamond+2\Bigl(\mu^m\,r\,\mat{\bD}(\vec U^{m+1,\ell+1}),~\mat{\bD}(\vec\chi)\Bigr)- \Bigl(P^{m+1,\ell+1},~\bG\cdot[r\vec\chi]\Bigr) \nn\\
&\hspace{1.5cm} -\,\gamma\,\big\langle(\vec X^{m}\cdot\vec e_1)\,\varkappa^{m+1,\ell+1},~\vec\nu^m\cdot\vec\chi\,|\vec X^m_\alpha|\big\rangle  = \Bigl(\rho^m\vec g,~\vec\chi\,r\Bigr)\qquad\forall\vec\chi\in\bU^m,\\[0.5em]
\label{eq:picard2}
&\hspace{0.1cm}\Bigl(\bG\cdot[r\vec U^{m+1,\ell+1}],~q\Bigr)=0\qquad\forall q\in \bP^m, \\[0.5em]
\label{eq:picard3}
&\frac{1}{\ttau}\big<\vec X^{m+1,\ell+1}-\vec X^m,~\zeta\,\vec f^{m+\frac{1}{2},\ell}\big>-\big<(\vec X^m\cdot\vec e_1)\vec U^{m+1,\ell+1}\cdot\vec\nu^m,~\zeta\,|\vec X^m_\alpha|\big>=0\qquad\forall\zeta\in V^h,\\[0.5em]
&\big<\varkappa^{m+1,\ell+1},~\vec\eta\cdot\vec f^{m+\frac{1}{2},\ell}\big> +\big<\vec\eta\cdot\vec e_1,~|\vec X^{m+1,\ell}_\alpha|\big> + \big<(\vec X^m\cdot\vec e_1)\vec X^{m+1,\ell+1}_\alpha,~\vec\eta_\alpha\,|\vec X^m_\alpha|^{-1}\big>=0\qquad\forall\vec\eta\in V_\partial^h,
\label{eq:picard4}
\end{align}
\end{subequations}
where $\vec f^{m+\frac{1}{2}, \ell}$ is a lagged approximation which follows \eqref{eq:averagenor} directly, except that $\vec X^{m+1}$ is replaced by $\vec X^{m+1, \ell}$. The equations in \eqref{eqn:picard} lead to a linear system which can be written in a form which is very similar to \eqref{eq:linsys}, and thus they can be solved via the introduced techniques. In fact, in the new linear system, we need to have $\mathcal{\vec N}_\Gamma^m$ replaced by $\vec N_\Gamma^m$, while on the right hand side, we have additional contribution from the second term in \eqref{eq:picard4}. We repeat the above iterations for $\ell = 0,\cdots$ until the following condition is satisfied
\begin{equation}
 \max_{1\leq j\leq J_\Gamma}|\vec X^{m+1,\ell +1}(\alpha_j) - \vec X^{m+1,\ell}(\alpha_j)|\leq {\rm tol},\label{eq:pictol}
\end{equation}
where ${\rm tol}$ is the chosen numerical tolerance.

In the following, we implement the introduced ALE methods by considering the experiments of a rising bubble and an oscillating droplet. Unless otherwise stated, we always choose $\vec U^0=\vec 0$ and ${\rm tol}=10^{-8}$ in \eqref{eq:pictol}. We also introduce the discrete quantities
\begin{equation}
{\alpha_{\rm min}}|_{t_m}=\min_{\sigma\in\mathscr{T}^m}\min_{\alpha\in \measuredangle{(\sigma)}},\qquad \Psi_e|_{t_m} := \frac{\max_{j=1}^{J_\Gamma}|\vec X^m(\alpha_j) - \vec X^m(\alpha_{j-1})|}{\min_{j=1}^{J_\Gamma}|\vec X^m(\alpha_j) - \vec X^m(\alpha_{j-1})|},\qquad v_\Delta|_{t_m}: = \frac{\vol(\vec X^m)}{\vol(\vec X^0)}-1,\nn
\end{equation}
where $\measuredangle{(\sigma)}$ is the set of all the angles of the simplex $\sigma$. Here $\alpha_{\rm min}$ and $\Psi_e$ measure the quality of the bulk mesh and interface mesh, respectively, and $v_\Delta$ is the relative volume loss of the inner phase. In practice, we observe that the interface mesh is well preserved, especially for the ALE methods that enjoy the property of equidistribution. Besides, the moving mesh approach described in \S\ref{sec:dale} in general works smoothly except when the interface exhibits strong deformations. To avoid the mesh deterioration, we regenerate the bulk mesh once the following condition is violated
\[\alpha_{\rm min}\geq \frac{\pi}{18}.\] 
After the remeshing, we then need to appropriately interpolate the fluid velocity and bulk mesh velocity to the new generated mesh so that the discrete ALE mappings in \eqref{eq:DALE} are well defined.

\subsection{The rising bubble}

We study the dynamics of a rising bubble, which was considered in the case of 2d in \cite{Hysing2009}, see also the generalization to 3d in \cite{BGN15stable},
and to the rotationally symmetric setting in \cite{GNZ23asy}. The physical parameters are given by  
\begin{itemize}
\item Case I:
\begin{equation}
\rho_+=1000, \quad\rho_- = 100, \quad\mu_+=10,\quad\mu_-=1,\quad \gamma = 24.5,\quad \vec g = (0, -0.98)^T.\nn
\end{equation}
\item Case II:
\begin{equation}
\rho_+=1000, \quad\rho_- = 1, \quad\mu_+=10,\quad\mu_-=0.1,\quad \gamma = 1.96,\quad \vec g = (0, -0.98)^T.\nn
\end{equation}
\end{itemize}
We follow the numerical setting from \cite{GNZ23asy} and use  $\mR = [0, \frac{1}{2}]\times[0,2]$ with $\partial_1\mR=[0,0.5]\times\{0,2\}$, $\partial_2\mR=\{0.5\}\times[0,2]$. The initial interface is given by $\Gamma(0)=\left\{\vec\mZ\in\mR:\;|\vec\mZ-(0, \frac{1}{2})^T|=\frac{1}{4}\right\}$. To monitor the state of the rising bubble, we introduce the discrete benchmark quantities 
\begin{equation}
\cira|_{t_m} := \frac{\pi^{\frac{1}{3}}[6\,M(\vec X^m)]^{\frac{2}{3}}}{A(\vec X^m)},\qquad V_c|_{t_m}:=\frac{2\pi\int_{\mR_-^m}(\vec U^m\cdot\vec e_2)r\,\drz}{\vol(\vec X^m)},\qquad 
z_c|_{t_m}:=\frac{2\pi\int_{\mR_-^m}(\vec\id\cdot\vec e_2)r\,\drz}{\vol(\vec X^m)},\nn
\end{equation}
where $\cira$ denotes the degree of sphericity, $V_c$ is the rise velocity, and $z_c$ is the centre of mass in the vertical direction.

\begin{table}[htp]
\centering
\def\temptablewidth{0.80\textwidth}
\vspace{0pt}
\caption{Benchmark quantities of the rising bubble in the case I, where $h=1/J_\Gamma$ with $h_0=2^{-4}$ and $\ttau_0=0.01$.}\label{tab:case1}
\renewcommand\arraystretch{1.25}
{\rule{\temptablewidth}{1pt}}
\begin{tabular}{c|ccc|ccc}
&\multicolumn{3}{c|}{n-$\Stab^h$ in \eqref{eqn:stabfd}} &\multicolumn{3}{c}{n-$\Equi^h$ in \eqref{eqn:meshfd}}\\ \hline 
$(h,\ttau)$  & $(h_0,\ttau_0)$   &$(\frac{h_0}{2}, \frac{\ttau_0}{4})$ &$(\frac{h_0}{4}, \frac{\ttau_0}{16})$ & $(h_0,\ttau_0)$   &$(\frac{h_0}{2}, \frac{\ttau_0}{4})$ &$(\frac{h_0}{4}, \frac{\ttau_0}{16})$  \\\hline  
$\cira_{\min}$ &0.9490 &0.9484  &0.9483 &0.9473  &0.9480 &0.9482 \\\hline 
$t_{\cira = \cira_{\min}}$ &3.0000 &3.0000  &3.0000 & 3.0000  &3.0000 &3.0000\\[0.4em]\hline 
$V_{c,\max}$ & 0.3689 &0.3664  &0.3657 & 0.3677 &0.3661 &0.3656  \\\hline 
$t_{_{V_c = V_{c,\max}}}$ & 0.9100 &0.9100  &0.9119 & 0.9200  &0.9150 &0.9125 \\\hline 
$z_c(t=3)$ &1.4925 &1.4890  &1.4882 &1.4892 &1.4883 &1.4880 \\ \hline 
$v_\Delta(t=3)$ &-8.80E-4 &-2.28E-4 &-5.74E-5 &-4.10E-4 &-1.07E-4 &-2.71E-5
 \end{tabular}
{\rule{\temptablewidth}{1pt}}
{\rule{\temptablewidth}{1pt}}
\begin{tabular}{c|ccc|ccc}
&\multicolumn{3}{c|}{c-$\Stab^h$ in \eqref{eqn:cstabfd}} &\multicolumn{3}{c}{c-$\Equi^h$ in \eqref{eqn:cmeshfd}}\\ \hline 
 $(h,\ttau)$ & $(h_0,\ttau_0)$   &$(\frac{h_0}{2}, \frac{\ttau_0}{4})$ &$(\frac{h_0}{4}, \frac{\ttau_0}{16})$ & $(h_0,\ttau_0)$   &$(\frac{h_0}{2}, \frac{\ttau_0}{4})$ &$(\frac{h_0}{4}, \frac{\ttau_0}{16})$   \\\hline  
$\cira_{\min}$ &0.9490 &0.9484  &0.9483 &0.9473   &0.9480 &0.9482 \\\hline 
$t_{\cira = \cira_{\min}}$ &3.0000 &3.0000  &3.0000 & 3.0000  &3.0000 &3.0000\\[0.4em]\hline 
$V_{c,\max}$ &0.3689 &0.3664  &0.3657 & 0.3677 &0.3661 &0.3656  \\\hline 
$t_{_{V_c = V_{c,\max}}}$ & 0.9100 &0.9100  &0.9119 & 0.9150 &0.9125 &0.9125  \\\hline 
$z_c(t=3)$ &1.4924 &1.4890  &1.4882 &1.4891 &1.4883 &1.4880\\ \hline 
$v_\Delta(t=3)$ &-8.80E-4 &-2.28E-4 &-5.74E-5 &-4.10E-4 &-1.07E-4 & -2.71E-5
   \end{tabular}
{\rule{\temptablewidth}{1pt}}
\end{table}

\medskip
\noindent
{\bf Example 1}: We first simulate the rising bubble in case I, using the introduced ALE methods from \S\ref{sec:nALEm} and \S\ref{sec:cALEm}. We employ three different computational meshes for each considered method, and the discrete benchmark quantities are reported in Table~\ref{tab:case1}. Based on the data, we can conclude that (i) the nonconservative and conservative ALE methods can produce very similar numerical results, despite the different treatments of the inertia terms in the two-phase Navier-Stokes equations; (ii) numerical convergence is observed as the mesh is refined, and the relative volume loss exhibits a second-order convergence. As the numerical results for the nonconservative and conservative methods are almost the same and graphically indistinguishable, in the following presentations we will only show the results from the nonconservative ALE methods.

\begin{figure}[t]
\centering
\includegraphics[width=0.9\textwidth]{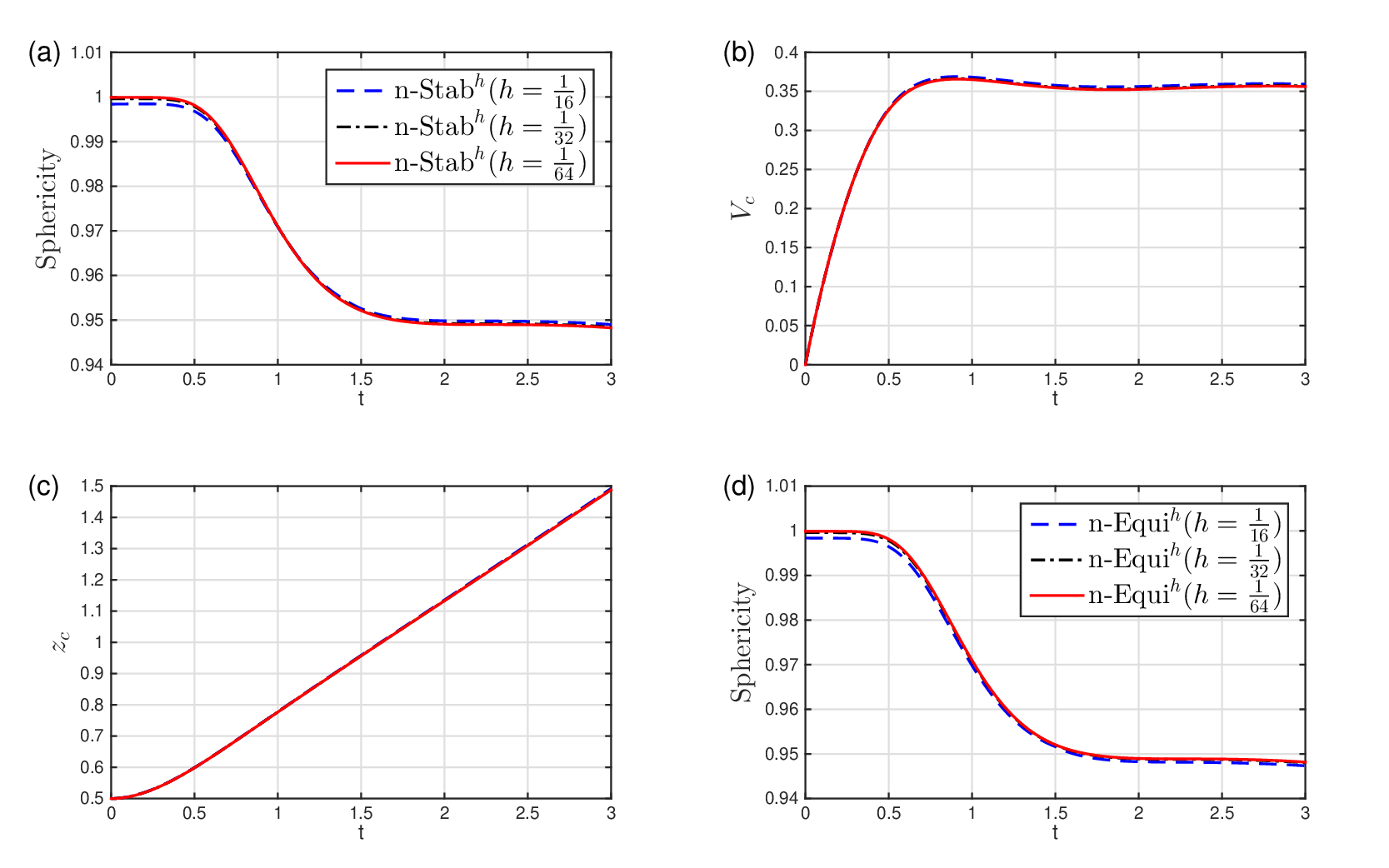}
\caption{The time history of the benchmark quantities, where (a),(b),(c) are from the method n-$\Stab^h$, and (d) is from the method n-$\Equi^h$.}
\label{Fig:B1Q}
\end{figure}

\begin{figure}[!htp]
\centering
\includegraphics[width=0.42\textwidth]{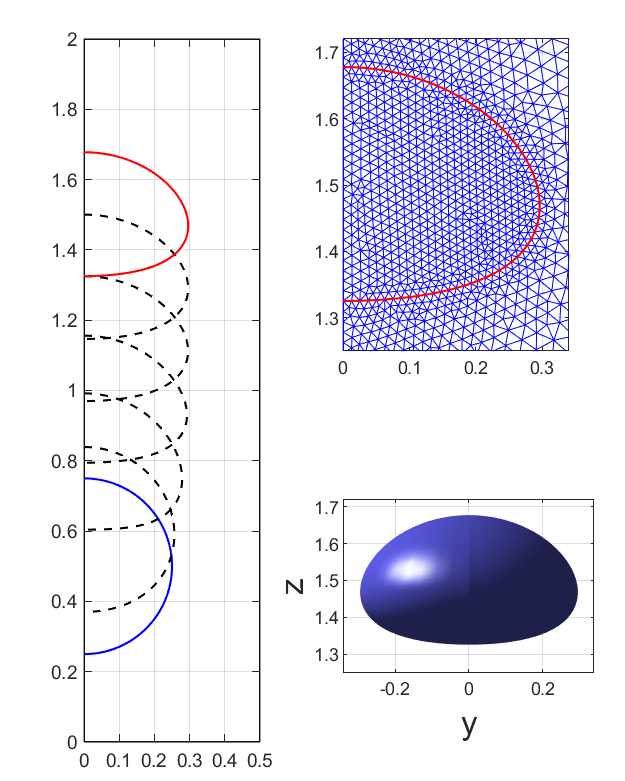}
\includegraphics[width=0.42\textwidth]{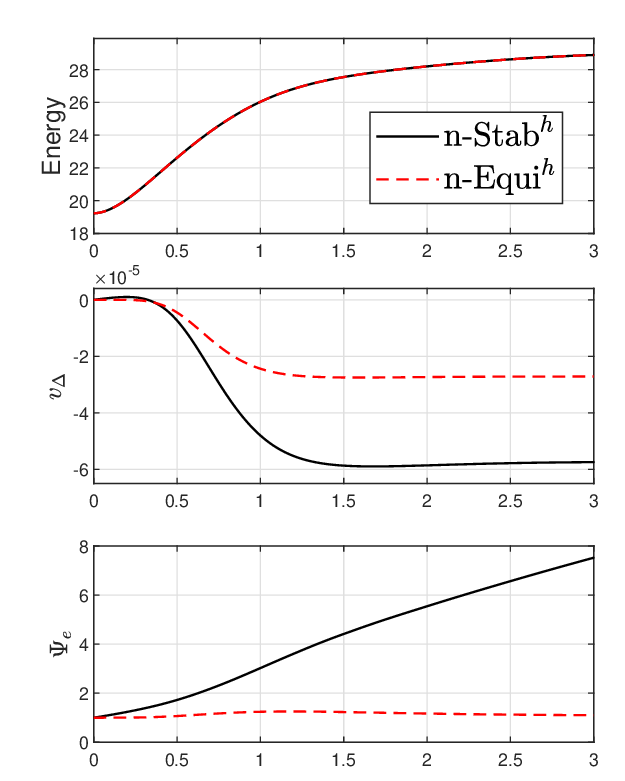}
\caption{Evolution of the rising bubble in case I. On the left we show the snapshots of the fluid interface at times $t=0,0.5,\cdots,3.0$ with visualisations of the computational mesh and the bubble at the final time. On the right are plots of the energy, the relative volume loss and the mesh quality indicator. Here $h=\frac{1}{64}$ and $\ttau=1.5625\times 10^{-4}$.}
\label{Fig:B1V}
\end{figure}

The benchmark quantities versus time are plotted in Fig.~\ref{Fig:B1Q}, which further verifies the numerical convergence. Visualisations of the fluid interfaces at several times are shown in Fig.~\ref{Fig:B1V}. In the same figure, we also plot the time history of the energy, the relative volume loss and the mesh quality indicator. Here we observe an excellent agreement of the energy from the methods n-$\Stab^h$ and n-$\Equi^h$, and the volume of the inner phase is nearly preserved.   Moreover, we find that the mesh indicator $\Psi_e$ gradually approaches 1 for the method n-$\Equi^h$, which implies the property of equidistribution. On the other hand, for the method n-$\Stab^h$ we observe that $\Psi_e$ gradually grows, without exceeding a value of 8. This implies that that the interface mesh quality for the method n-$\Stab^h$ is also preserved well.

\begin{figure}[!htp]
\centering
\includegraphics[width=0.90\textwidth]{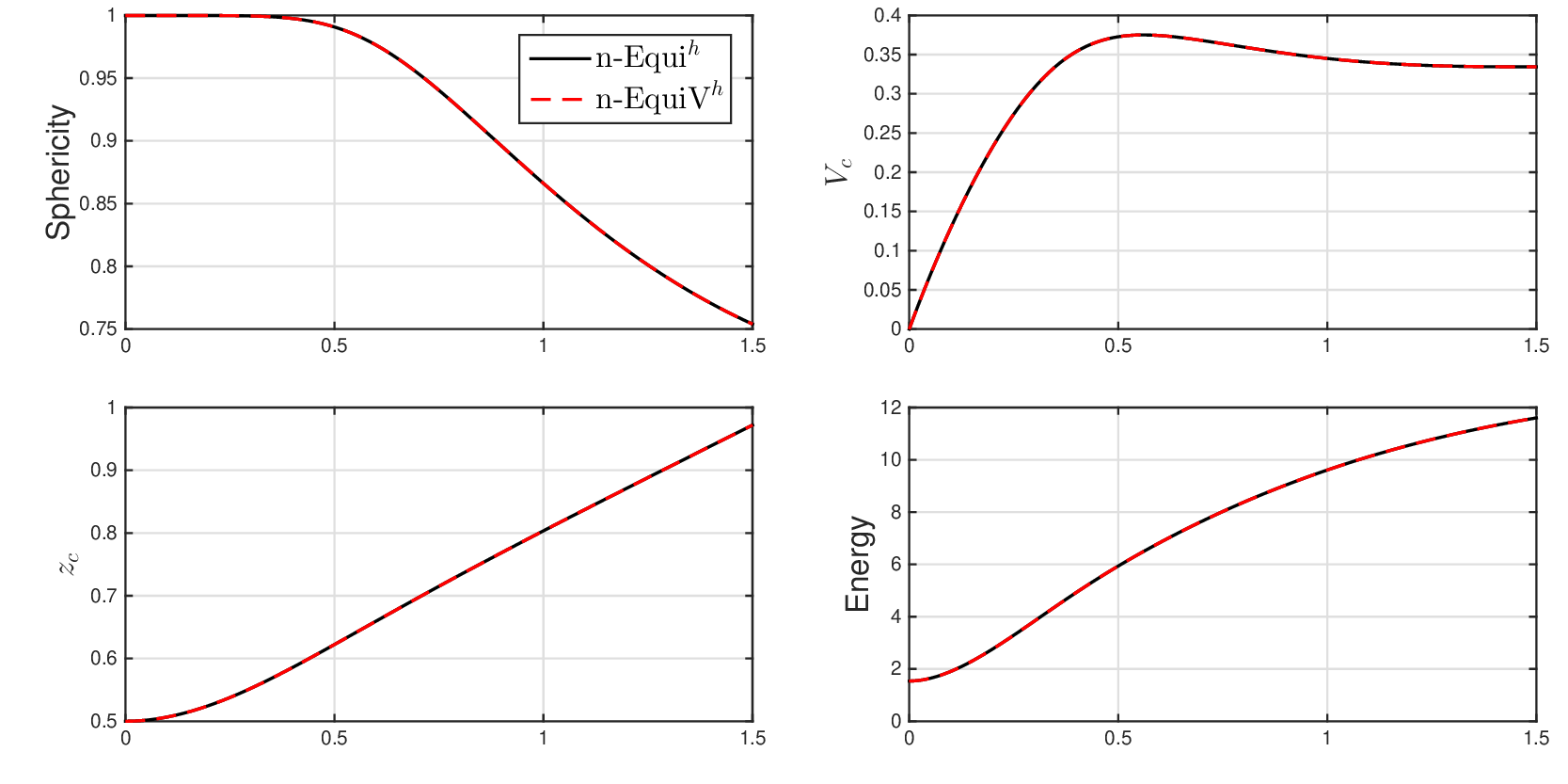}
\caption{The time history of the discrete benchmark quantities for the methods n-$\Equi^h$ and n-$\EquiV^h$, where $h=\frac{1}{64}$ and $\ttau=1.5625\times 10^{-4}$. }
\label{Fig:B2Q}
\end{figure}

\begin{figure}[!htp]
\centering
\includegraphics[width=0.90\textwidth]{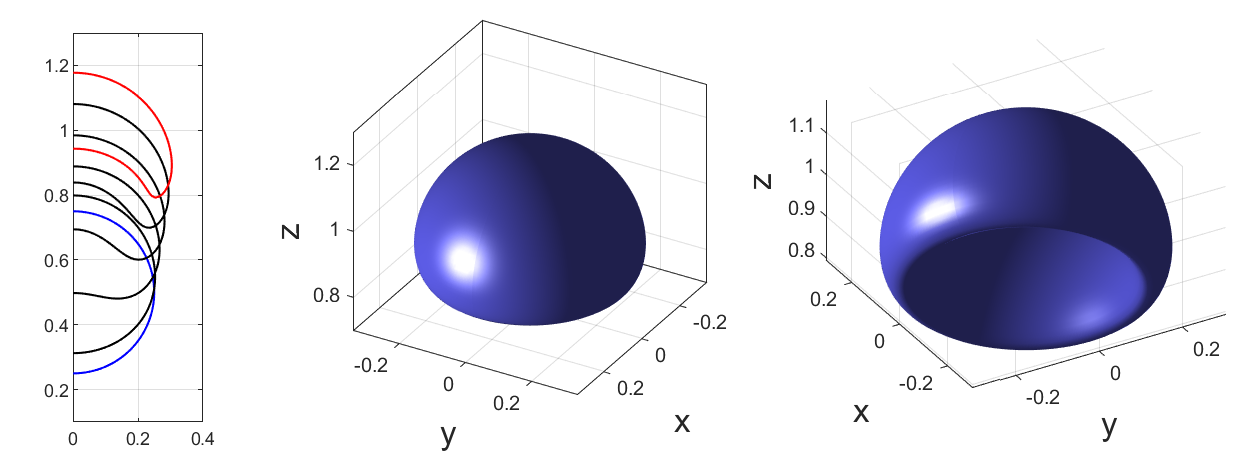}
\caption{Evolution of the rising bubble in case II using the method n-$\EquiV^h$, where we show the generating curves at times $t=0,0.3,0.6,0.9,1.2,1.5$ and the axisymmetric interfaces $\mS(t)$ at $t=1.5$ with views from different directions. Here $h=\frac{1}{64}$, $\ttau=1.5625\times 10^{-4}$.}
\label{Fig:B2V}
\end{figure}

\begin{figure}[!htp]
\centering
\includegraphics[width=0.90\textwidth]{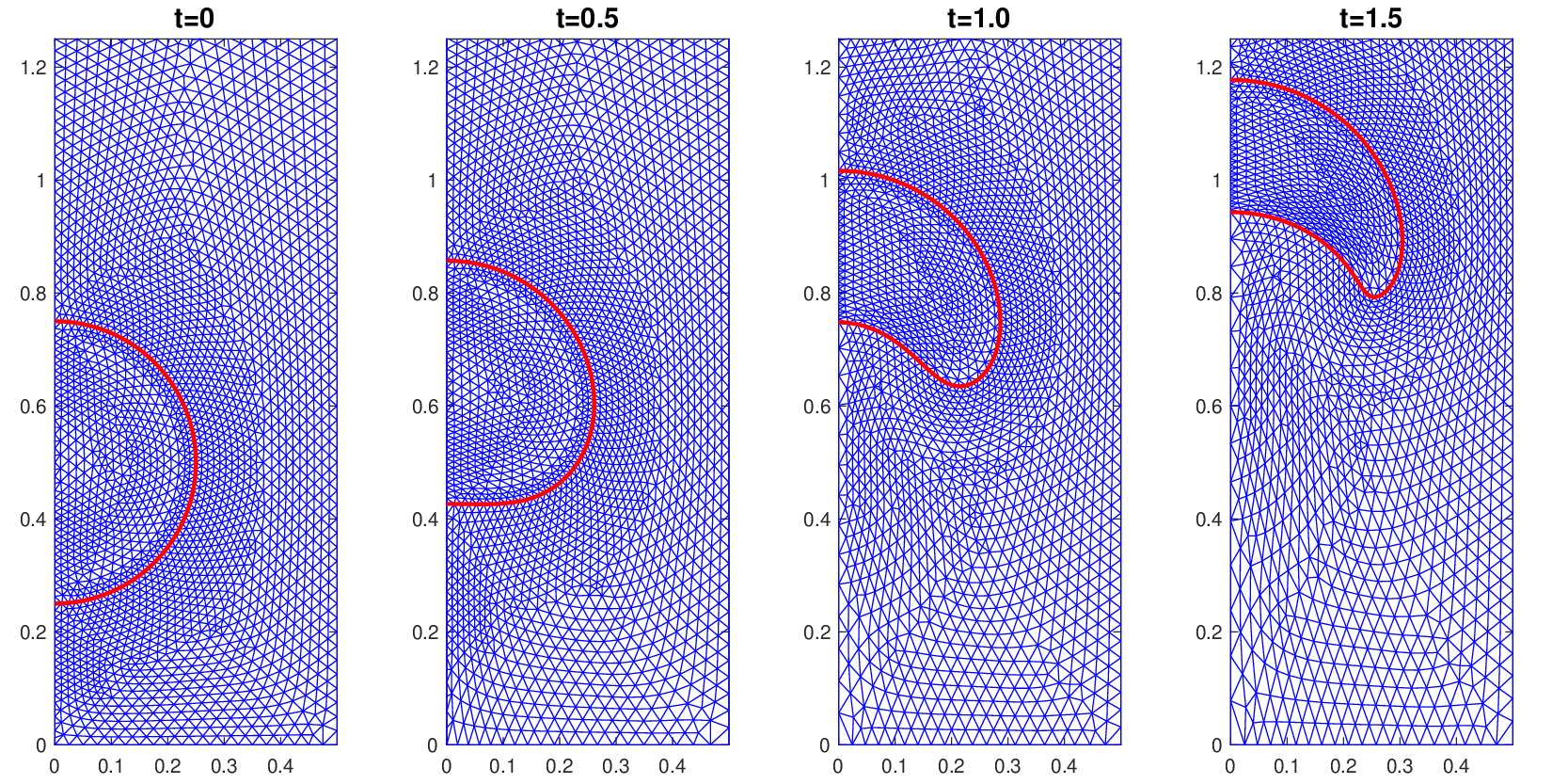}
\caption{The computational meshes at several times $t=0, 0.5, 1.0, 1.5$ for the experiment in Fig.~\ref{Fig:B2V}.}
\label{Fig:B2Mesh}
\end{figure}

\medskip
\noindent
{\bf Example 2}: We next study the dynamics of the rising bubble in case II by using the methods introduced in \S\ref{sec:nALEm}, and their volume-preserving variants from \S\ref{sec:vp}. Once again we consider three different sets of discretization parameters for each method, and the numerical results are reported in Table~\ref{tab:case2}. Here we observe that these methods can produce very similar numerical results. In particular, the volume is exactly preserved for the methods n-$\StabV^h$ and n-$\EquiV^h$, which numerically verifies \eqref{eq:DVC} and \eqref{eq:DVCC}. The benchmark quantities and the evolving fluid interface are shown in Fig.~\ref{Fig:B2Q} and Fig.~\ref{Fig:B2V}, respectively.  Several snapshots of the corresponding computational mesh are shown in Fig.~\ref{Fig:B2Mesh}.

\begin{table}[!htp]
\centering
\def\temptablewidth{0.8\textwidth}
\vspace{0pt}
\caption{Benchmark quantities of the rising bubble in case II, where $h=1/J_\Gamma$ with $h_0=2^{-4}$ and $\ttau_0=0.01$.}
\label{tab:case2}
\renewcommand\arraystretch{1.25}
{\rule{\temptablewidth}{1pt}}
\begin{tabular}{c|ccc|ccc}
&\multicolumn{3}{c|}{n-$\Stab^h$} &\multicolumn{3}{c}{n-$\Equi^h$}\\ \hline 
$(h,\ttau)$  & $(h_0,\ttau_0)$   &$(\frac{h_0}{2}, \frac{\ttau_0}{4})$ &$(\frac{h_0}{4}, \frac{\ttau_0}{16})$ & $(h_0,\ttau_0)$   &$(\frac{h_0}{2}, \frac{\ttau_0}{4})$ &$(\frac{h_0}{4}, \frac{\ttau_0}{16})$  \\\hline  
$\cira_{\min}$ &0.7504 &0.7538  &0.7547 & 0.7411 & 0.7511  & 0.7540  \\\hline 
$t_{\cira = \cira_{\min}}$ &1.5000 &1.5000  &1.5000 & 1.5000  &1.5000 &1.5000\\[0.4em]\hline 
$V_{c,\max}$ &0.3787 &0.3758  &0.3750  & 0.3780 &0.3756 &0.3750   \\\hline 
$t_{_{V_c = V_{c,\max}}}$ & 0.5600 &0.5600  &0.5600 & 0.5600  &0.5600 &0.5600  \\\hline 
$z_c(t=1.5)$ &0.9723 &0.9722  &0.9721 &0.9710 &0.9719 &0.9720  \\ \hline 
$v_\Delta(t=1.5)$ &-2.95E-3 &-6.89E-4 &-1.77E-4 &-6.96E-4 &-2.04E-4 & -5.14E-5 
 \end{tabular}
{\rule{\temptablewidth}{1pt}}
{\rule{\temptablewidth}{1pt}}
\begin{tabular}{c|ccc|ccc}
&\multicolumn{3}{c|}{n-$\StabV^h$} &\multicolumn{3}{c}{n-$\EquiV^h$}\\ \hline 
$(h,\ttau)$  & $(h_0,\ttau_0)$   &$(\frac{h_0}{2}, \frac{\ttau_0}{4})$ &$(\frac{h_0}{4}, \frac{\ttau_0}{16})$ & $(h_0,\ttau_0)$   &$(\frac{h_0}{2}, \frac{\ttau_0}{4})$ &$(\frac{h_0}{4}, \frac{\ttau_0}{16})$  \\\hline  
$\cira_{\min}$ &0.7514 &0.7541 &0.7548 &0.7415 &0.7514  &0.7540  \\\hline 
$t_{\cira = \cira_{\min}}$ &1.5000 &1.5000  &1.5000 & 1.5000  &1.5000 &1.5000\\[0.4em]\hline 
$V_{c,\max}$  &0.3792 &0.3759 &0.3751  &0.3786 &0.3758  &0.3750 \\\hline 
$t_{_{V_c = V_{c,\max}}}$ & 0.5600 &0.5600  &0.5600 & 0.5700  &0.5600 &0.5600  \\\hline 
$y_c(t=1.5)$  &0.9712 &0.9719 &0.9720   &0.9707 &0.9714  &0.9720 \\ \hline 
$v_\Delta(t=1.5)$ &5.22E-9 &3.16E-9 &1.09E-9 &-4.75E-10 &5.91E-10 &3.59E-10 
 \end{tabular}
{\rule{\temptablewidth}{1pt}}
\end{table}

\vspace{0.4cm}

\subsection{The oscillating droplet}
In this subsection, we numerically study the oscillation of a levitated drop which is surrounded by a fluid of low density. Inspired by the work in \cite{aalilija20}, we consider an axisymmetric perturbation of a spherical equilibrium. In particular, the generating curve of the initial drop is given by   
\begin{equation}
\left\{\begin{array}{ll}r(\theta, 0) &= R_0\left[1 + \varepsilon_{n,0} P_n(\cos\theta)- \frac{1}{2n+1}\varepsilon_{n,0}^2\right]\cos(\theta-\frac{\pi}{2}),\\
z(\theta, 0) &= R_0\left[1 + \varepsilon_{n,0} P_n(\cos\theta)- \frac{1}{2n+1}\varepsilon_{n,0}^2\right]\sin(\theta-\frac{\pi}{2}) + 1.0,
\end{array}\right.\quad \theta\in[0,\pi],\quad n\geq 2,\label{eq:IniG}
\end{equation}
where $R_0$ is the radius, $\varepsilon_{n,0}$ is the magnitude of the perturbation, and $P_n(x)$ are Legendre polynomials. For example, $P_2(x)=\frac{1}{2}(3x^2-1)$ and $P_5(x)=\frac{1}{8}(63x^5 -70x^3 + 15x)$.  
Then on recalling the analytical asymptotic solution \cite[(15b) and (38)]{aalilija20}, we note that the dynamic generating curve can be approximated by
\begin{equation}
\left\{\begin{array}{ll}r(\theta, t) &= R_0\left[1 + \varepsilon_n(t)P_n(\cos\theta)- \frac{1}{2n+1}\varepsilon_n^2(t)\right]\cos(\theta-\frac{\pi}{2}),\\
z(\theta, t) &= R_0\left[1 + \varepsilon_n(t) P_n(\cos\theta)- \frac{1}{2n+1}\varepsilon_n^2(t)\right]\sin(\theta-\frac{\pi}{2}) + 1.0,
\end{array}\right.\quad \theta\in[0,\pi],\quad t\geq 0,\quad n\geq 2,\label{eq:otf}
\end{equation}
where $\varepsilon_n(t)$ is given by
\begin{equation}
\varepsilon_n(t)\approx \varepsilon_{n,0}\exp(-\lambda_{n} t)\cos(\omega_{n} t)\quad\mbox{with}\quad \omega_n = \sqrt{\omega_{n,0}^2-\lambda_n^2},\nn
\end{equation}
and
\begin{equation}
\omega_{n,0}=\sqrt{\frac{n(n-1)(n+2)\gamma}{\rho_-R_0^3}},
\qquad \lambda_n = \frac{(2n+1)(n-1)\mu_-}{\rho_-R_0^2}.\nn
\end{equation}
We further introduce the radius of the droplet as $R(\theta, t) = R_0\left[1 + \varepsilon_n(t) P_n(\cos\theta)- \frac{1}{2n+1}\varepsilon_n^2(t)\right]$. For this experiment we use the computational domain $\mR=[0,0.6]\times[0,2]$ with $\partial_1\mR =[0,0.6]\times\{0,2\}$ and $\partial_2\mR=\{0.6\}\times[0,2]$, and  choose 
\begin{equation}
\rho_+=1, \quad\rho_- = 1000, \quad\mu_+=0.01,\quad\mu_-=2,\quad \gamma = 40,\quad \vec g = \vec 0,\quad R_0=0.3.\nonumber
\end{equation}

\begin{figure}[!htp]
\centering
\includegraphics[width=0.9\textwidth]{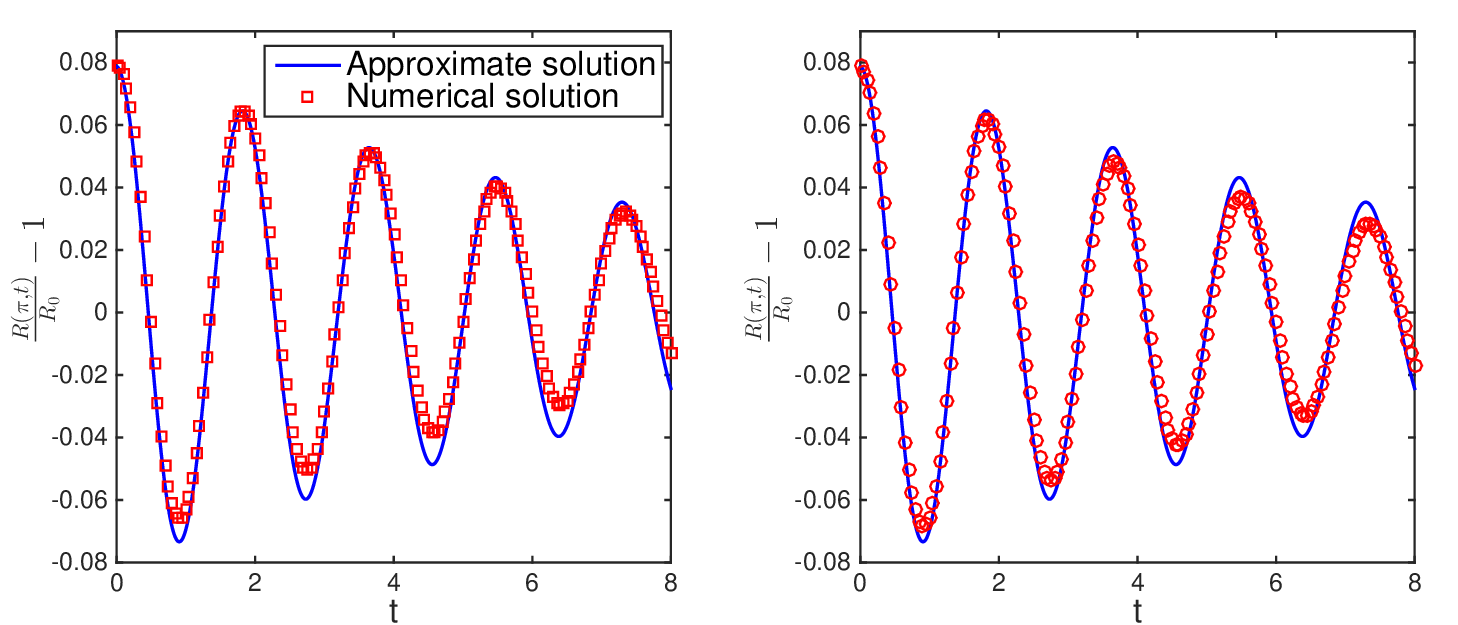}
\caption{[$n=2, \varepsilon_{n,0}=0.08$] The displacement for the upper point of the generating curve on the z-axis, where the numerical results are obtained by using the methods n-$\StabV^h$ (left panel) and n-$\EquiV^h$ (right panel) with $\ttau=10^{-3}$, $K=1348$, $J_\mR = 2574$ and  $J_\Gamma =64$. }
\label{Fig:Osci2}
\end{figure}

\medskip
\noindent
{\bf Example 3}:
We first focus on the $2$-mode perturbation of the droplet with $\varepsilon_{2,0} = 0.08$ in \eqref{eq:IniG} for the initial interface. The obtained numerical results by the method n-$\EquiV^h$ and n-$\StabV^h$ are compared with the approximate solution in \eqref{eq:otf} with $n=2$. As shown in Fig.~\ref{Fig:Osci2},  we observe an excellent agreement between the numerical solution and the approximate solution for both the introduced methods n-$\StabV^h$ and n-$\EquiV^h$. 

\begin{figure}[!htp]
\centering
\includegraphics[width=0.85\textwidth]{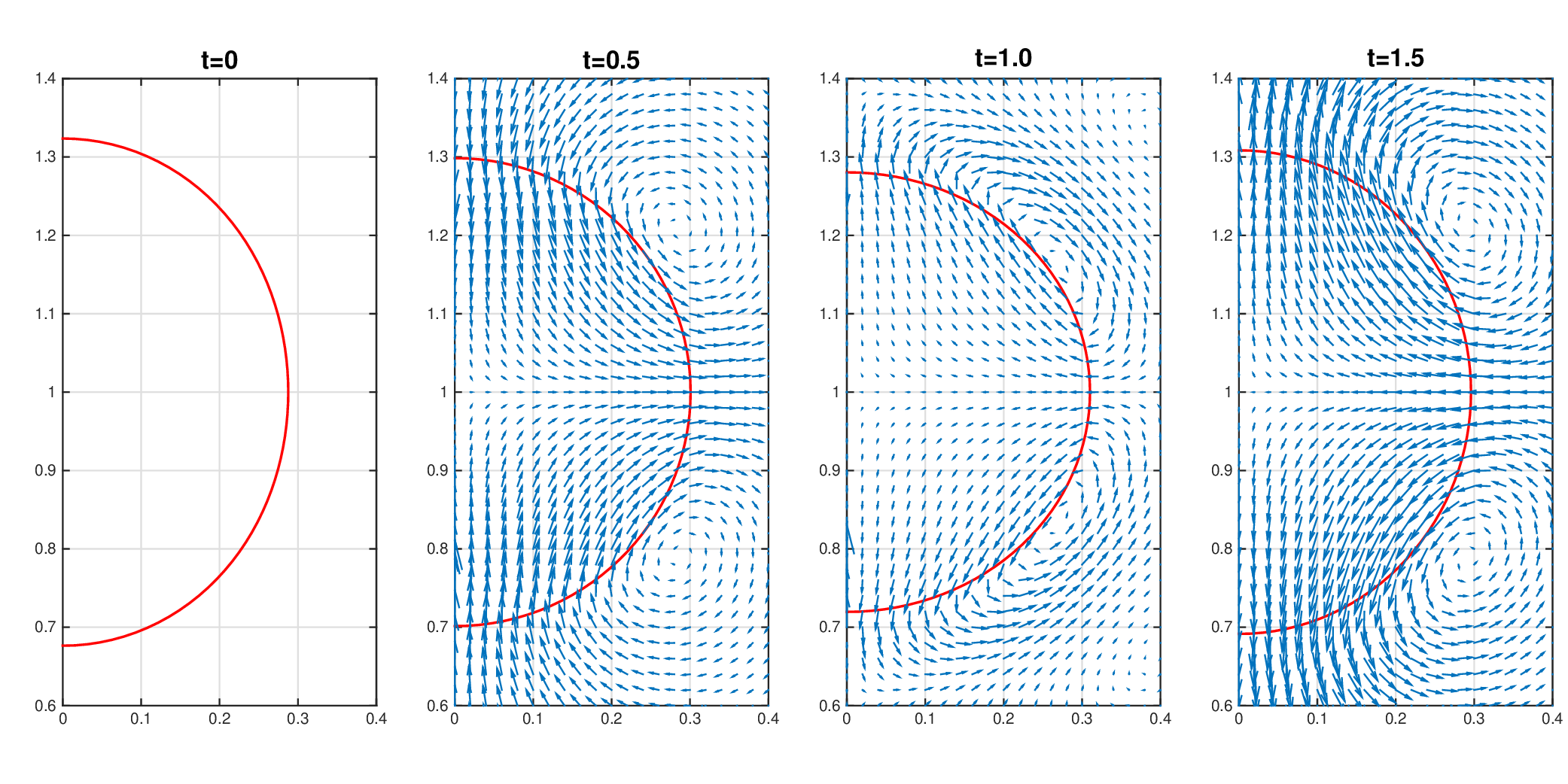}
\caption{Snapshots of the fluid interface and the velocity fields at several times for the experiment in Fig.~\ref{Fig:Osci2}.}
\label{Fig:Osci2v}
\end{figure}

\begin{figure}[!htp]
\centering
\includegraphics[width=0.85\textwidth]{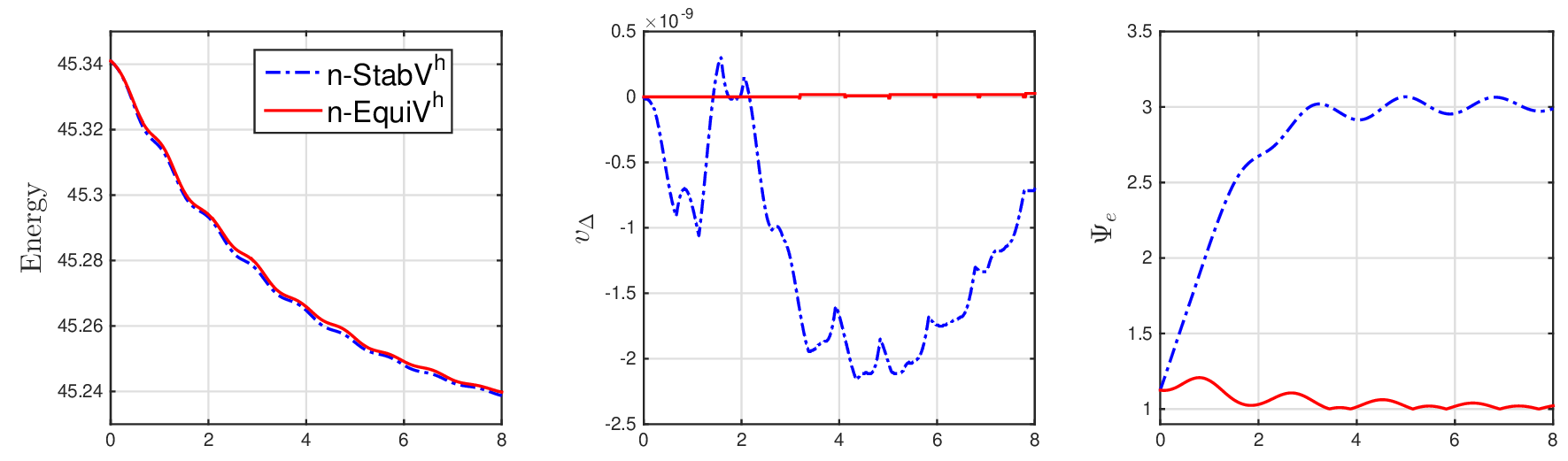}
\caption{The time history of the discrete quantities for the experiment in Fig.~\ref{Fig:Osci2}.}
\label{Fig:Osci2q}
\end{figure}

Snapshots of the interface and velocity fields are visualized in Fig.~\ref{Fig:Osci2v}. At time $t=0.5$, we observe that a counterclockwise vortex is generated in the upper half region, and a clockwise vortex is generated at the bottom. This implies that the droplet is spreading in the horizontal direction. While at $t=1.0$ and $t=1.5$, the droplet is spreading in the vertical direction. The benchmark quantities are plotted in Fig.~\ref{Fig:Osci2q}, which further verifies the good properties of the introduced methods.  

\begin{figure}[!htp]
\centering
\includegraphics[width=0.85\textwidth]{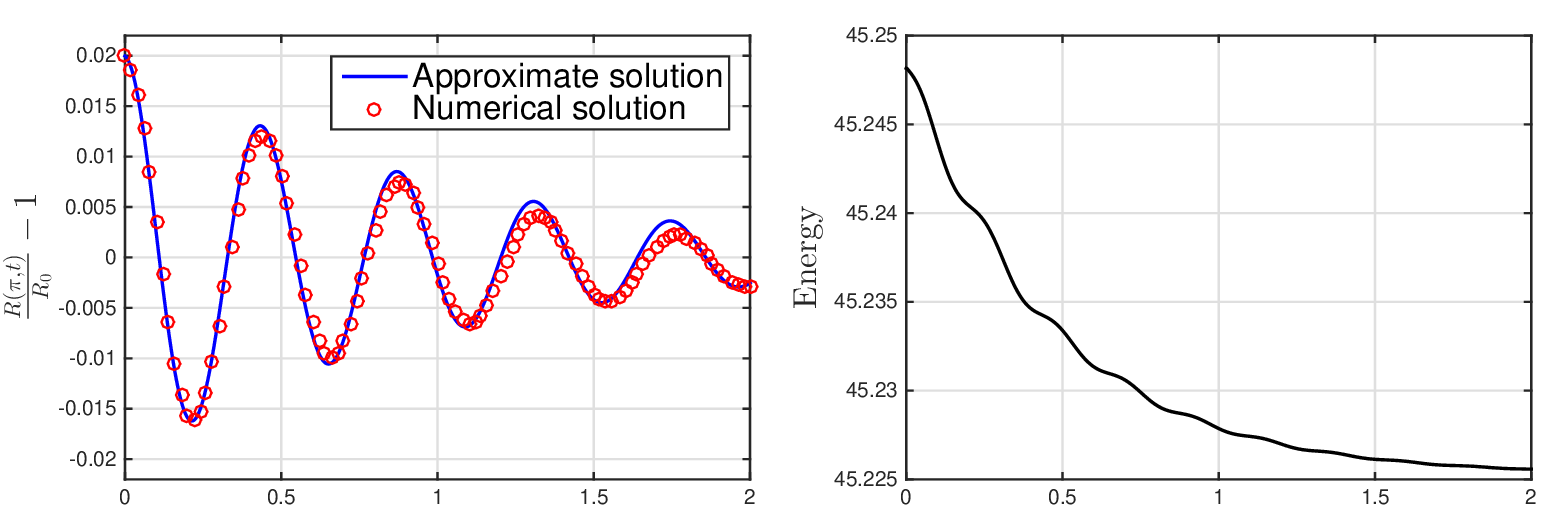}
\caption{[$n=5, \varepsilon_{n,0}=0.02$] The displacement for the upper point of the generating curve on the z-axis, where the numerical results are obtained by using the method n-$\EquiV^h$ (right panel) with $\ttau=5\times 10^{-4}$, $K=2993$, $J_\mR = 5804$ and  $J_\Gamma =64$. }
\label{Fig:5mode}
\end{figure}

\medskip
\noindent
{\bf Example 4}:  To further test the accuracy of our introduced ALE methods, we next consider the case of a $5$-mode perturbation. The same computational parameters are used except that the initial generating curve is given by \eqref{eq:IniG} with $n=5$, and $\varepsilon_{n,0}=0.02$. The numerical results for $n=5$ are shown in Fig. \ref{Fig:5mode}, where we observe an excellent agreement as well.

\section{Conclusions}\label{sec:con}
We proposed and analyzed a variety of ALE finite element approximations for the axisymmetric two-phase Navier-Stokes flow in both the conservative and nonconservative form. The introduced methods were shown to satisfy either unconditional stability or the equidistribution property, relying on two different approximations of the mean curvature of the axisymmetric interface. With the help of time-weighted approximations of the interface normals, we further adapted the introduced methods to achieve exact volume preservation on the discrete level. Numerical examples for a rising bubble and oscillating droplet were provided to examine the performance of these introduced methods. We observed that these methods can produce very accurate results, and that the volume preservation and energy stability are satisfied well on the discrete level. In future, we will consider the application of the introduced methods to more complex problems in two-phase flow.

\section*{Acknowledgements}
The last author would like to acknowledge the support from the Alexander von Humboldt Foundation during his stay in University of Regensburg. 
\appendix 
%\numberwithin{equation}{section}
\setcounter{equation}{0}

\section{Differential calculus}
Let $\varphi:\mR\times[0,T]\to\bR$ be a scalar field.  Applying the Reynolds transport theorem on $\mR_\pm(t)$ in terms of the ALE moving frame, it is not difficult to show that
\begin{align}
\ddt\int_{\mR_\pm(t)}\varphi\,r\,\drz&= \ddt\int_{\mR_\pm(t)}\varphi\,[\vec x\cdot\vec e_1]\,\drz \nn\\
&=\int_{\mR_\pm(t)}\bigl(\partial_t^\circ\varphi\,r + \varphi\,[\vec w\cdot\vec e_1]\bigr)\drz +  \int_{\mR_\pm(t)}\varphi\,r\,\bG\cdot\vec w\,\drz\nn\\
&=\int_{\mR_\pm(t)}\bigl(\partial_t^\circ\varphi\,r + \varphi\nabla\cdot[r\,\vec w]\bigr)\drz,\label{eq:ddtALEframe} 
\end{align}
where $\partial_t^\circ$ is defined in \eqref{eq:meshD} as the derivative with respect to the ALE moving reference, and $\vec w$ is the mesh velocity defined in \eqref{eq:ALEmeshv}. Applying $\varphi = \rho_\pm\vec u\cdot\vec\chi$ to \eqref{eq:ddtALEframe} and combining the two equations yields that
\begin{equation}
\ddt\Bigl(\rho\,\vec u,~\vec\chi\,r\Bigr) = \Bigl(\rho\,\partial_t^\circ\vec u,~\vec\chi\,r\Bigr) + \Bigl(\rho\,\partial_t^\circ\vec\chi,~\vec u\,r\Bigr) + \Bigl(\rho\,\vec u\cdot\vec\chi,~\nabla\cdot[r\,\vec w]\Bigr),\label{eq:ddtchiu}
\end{equation}
where $(\cdot,\cdot)$ is the $L^2$--inner product over $\mR$. This immediately implies that
\begin{equation}
\frac{1}{2}\ddt\Bigl(\rho\,\vec u,~\vec u\,r\Bigr) = \Bigl(\rho\,\partial_t^\circ\vec u,~\vec u\,r\Bigr) + \frac{1}{2}\Bigl(\rho\,|\vec u|^2,~\nabla\cdot[r\,\vec w]\Bigr).\label{eq:ddtK}
\end{equation}

\bibliographystyle{model1b-num-names}
\bibliography{bib}
\end{document}